\documentclass[a4paper,10pt]{book}
\usepackage[cp1251]{inputenc}
\usepackage{epigraph}

 \usepackage{amsmath, amsfonts, amssymb}
 \usepackage{amscd}
 \usepackage[russian]{babel}
 \usepackage{makeidx}
 \usepackage{calrsfs}
 \usepackage{graphicx}
 \usepackage{epsfig}
 \usepackage{curves, eepic, epic}
 \usepackage{wrapfig}
 \usepackage{bbm}
 \input diagxy
 \makeindex

 \DeclareMathAlphabet{\mathbbt}{U}{bbold}{m}{n}

\begin{document}

 \makeatletter

 \newcommand{\Author}[2]{\begin{center}\textbf{\large #1} \end{center} \medskip
                    \renewcommand{\@evenhead}{\raisebox{1mm}[\headheight][0pt]%
                    {\vbox{\hbox to\textwidth{\thepage \hfill\strut {\small #2}\hfill}\hrule}}} }

 \newcommand{\shorttitle}[1]{\renewcommand{\@oddhead}{\raisebox{1mm}[\headheight][0pt]%
                    {\vbox{\hbox to\textwidth{\strut \hfill{\small #1}\hfill\thepage}\hrule}}} }
 \headsep=2mm

 \renewcommand{\section}{\@startsection{section}{1}{\parindent}
                        {3ex plus 1ex minus .2ex}{2ex plus .2ex}{\bf\centering}}
 \renewcommand{\l@section}[2]{\footnotesize\leftskip0pt\par\noindent\hangindent27pt\hangafter=1{\qquad#1} \dotfill~~#2}
 \renewcommand{\l@part}[2]{\normalsize\leftskip0pt\par\smallskip\noindent\hangindent7pt\hangafter=1{\bf #1} \normalsize\dotfill~~#2}

 \newcommand{\Adress}[1]{\par\bigskip\baselineskip=11pt\hangindent17pt\hangafter=0\noindent{\footnotesize#1}\par\normalsize}
 \newcommand{\UDC}[1]{\begingroup\newpage\thispagestyle{empty}\begin{flushleft}УДК #1\end{flushleft}}
 \newcommand{\Abstract}[1]{\hangindent17pt\hangafter=0\noindent{\footnotesize#1}\bigskip\par\medskip}

 \makeatother

 \newcommand{\Title}[1]{\begin{center}\large\uppercase{#1}\end{center}\par}

 \def\subsec#1{\smallskip\textbf{#1}}
 \newcommand{\bib}[2]{{\leftskip-10pt\baselineskip=11pt\footnotesize\item{}\textsl{#1}~#2\par}}

\newcommand{\Proclaim}[1]{\smallskip{\bf#1}\sl}
\newcommand{\proclaim}[1]{{\bf#1}\sl}
\newcommand{\Theorem}[1]{\smallskip\par\textbf{Теорема#1.}\sl}
\newcommand{\Sledstvie}[1]{\smallskip\par\textbf{Следствие#1.}\sl}
\newcommand{\Lemma}[1]{\smallskip\par\textbf{Лемма#1.}\sl}
\newcommand{\Predl}[1]{\smallskip\textbf{Предложение#1.}\sl}
\newcommand{\Utv}[1]{\smallskip\textbf{Утверждение#1.}\sl}
\newcommand{\Endproc}{\rm}
\newcommand{\Def}[1]{\smallskip\par{\sc Определение#1.}}
\newcommand{\Zam}[1]{\smallskip\par{\sc Замечание#1.}}
\newcommand{\Primer}[1]{\smallskip\par{\sc Пример#1}}
\newcommand{\Zad}[1]{{\sc Задача#1.}}
\newcommand{\Lema}[1]{\smallskip\textbf{{Lemma#1}.\/~}\sl}
\newcommand{\Teorem}[1]{\smallskip\textbf{{Theorem#1}.\/ }\sl}
\newcommand{\Corollary}[1]{\smallskip\textbf{{Corollary#1}.\/~}\sl}
\newcommand{\Proposition}[1]{\smallskip\textbf{{Proposition#1}.\/~}\sl}
\newcommand{\Definition}[1]{\smallskip{\sc Definition#1.}}
\newcommand{\Remark}[1]{\smallskip{\sc Remark#1.}~}

\def\beginproof{\par\mbox{$\vartriangleleft$}}
\def\endproof{\text{$\vartriangleright$}}

\renewcommand{\le}{\leqslant}
\renewcommand{\ge}{\geqslant}
\renewcommand{\to}{\rightarrow}

\def\osum{\mathop{o\text{\/-}\!\sum}}  
\def\bosum{bo\text{-}\!\sum}
\def\olim{\mathop{o\text{-}{\fam0 lim}}} 
\def\bolim{\mathop{bo\text{-}\!\lim}}
\def\rlim{\mathop{r\text{-}{\fam0 lim}}}

\def\card{\mathop{\fam0 card}}
\def\co{\mathop{\fam0 co}}
\def\Dom{\mathop{\fam0 Dom}}
\def\dom{\mathop{\fam0 dom}}
\def\End{\mathop{\fam0 End}}
\def\Fin{\mathop{\fam0 fin}\nolimits}
\def\fin{\mathop{\fam0 fin}\nolimits}
\def\Id{\mathop{\fam0 Id}}
\def\Im{\mathop{\fam0 Im}}
\def\ind{\mathop{\fam0 ind}}
\def\Ker{\mathop{\fam0 Ker}}
\def\loc{\mathop{\fam0 loc}}
\def\const{\mathop{\fam0 const}}
\def\Exp{\mathop{\fam0 Exp}}
\def\grad{\mathop{\fam0 grad}}
\def\proj{\mathop{\fam0 proj}}
\def\im{\mathop{\fam0 im}}
\def\Lip{\mathop{\fam0 Lip}}
\def\mod{\mathop{\fam0 mod}}
\def\On{\mathop{\fam0 On}}
\def\Orth{\mathop{\fam0 Orth}}
\def\rank{\mathop{\fam0 rank}}
\def\sa{\mathop{\fam0 sa}}
\def\SC{\mathop{\fam0 SC}}
\def\sign{\mathop{\fam0 sign}}
\def\St{\mathop{\fam0 St}}
\def\supp{\mathop{\fam0 supp}}
\def\Re{\mathop{\fam0 Re}}

\def\Cal#1{\mathcal{#1}}   
\def\goth#1{\mathfrak{#1}} 
\def\bbold#1{{\mathbbt #1}} 

\newcommand{\shortpage}{\enlargethispage{-\baselineskip}}
 \newcommand{\Signatura}[3]{\begin{picture}(0,0)\put(#2,#3){\makebox{{#1}}}\end{picture}}

 \newpage

\Title{
Generalized Young and Cauchy--Bunyakowsky  Inequalities with Applications:\\
A SURVEY.
}

\noindent\rule{0.4\linewidth}{0.3pt}

\noindent Originally published  in the book:\\

"Advances in Modern Analysis".\\

Mathematical Forum Series.\\

Editors: Yu.F.\,Korobeinik,  A.G.\,Kusraev.\\

Vladikavkaz: South Mathematical Institute of the Vladikavkaz Scientific Center
of the Russian Academy of Sciences and the Goverment of  Republic of North Ossetia--Alania .\\

2009.\  P. 221--266.

\shorttitle{Generalized Inequalities with Applications.}

\Author{Author: S.~M.~Sitnik}{Sitnik~S.~M.}

\addcontentsline{toc}{contpart}{\protect  Sitnik~S.~M.
{\rm Generalized Inequalities with Applications}}
\baselineskip=1.05\baselineskip

\begin{center}
Chair on   Mathematics,\\
Department of Radioengineering,\\
Voronezh Militia Institute
of the Ministry of Internal Affairs of Russia,\\
Voronezh, Russia.\\

e-mail: pochtasms@gmail.com\\

postal: P.O.Box 12, Voronezh-5, \        Voronezh, \ 394005, \ Russia.
\end{center}

\medskip

In this survey we consider generalizations for Young and Cauchy--Bunyakowsky  inequalities with applications.

\newpage

\begin{center}
INTRODUCTION.
\end{center}

\vskip 0.3cm
In this survey we consider generalizations for Young and Cauchy--Bunyakowsky  inequalities with applications.

In section 2 it is noticed that in fact there two different Young inequalities (for two numbers) but not one. Then there is a problem to compare different versions.   In theorem 2 three different cases are considered and three numerical examples for all cases added.

In section 3 generalizations are considered on  Young inequality for  integrals and for Legendre transforms, applications on improvement of Rogers--H\"older--Riesz inequality. We mention the \textit{transfere principle} suggested by the author and proved recently by A.G.\,Kusraev how to transfere inequalities for numbers to the case of uniformly complete vector lattices and operators in them.

In section 4 different kinds of means are  considered: power, Rado, Gini, Lemer, quasi--arithmetic, iterated and AGM and also different known procedures for constructing new means from  known ones. In particular in theorem 4 simplified version of Tibor Rado results on comparison of  power and Rado means is presented. After that abstract means are introduced and in theorem 5   \textit{all means} for two numbers are characterized completely.   The generalized entropy is introduced and   completely characterized in theorem 6.

In section 5 discrete version of  Cauchy--Bunyakowsky  inequality is studied and its generalizations: Milne, Callebaut and Daykin--Eliezer--Carlitz. The author suggest to express them in terms of means --- it is much more clear and give a lot of new non--trivial examples.

In section 6 integral version of  Cauchy--Bunyakowsky  inequality is studied with its improvements. In theorem
8 an integral version is proved of  sufficiency of Daykin--Eliezer--Carlitz result. Different special cases are considered for power, Rado, $\min-\max$ and AGM means, the last one is quiet exotic. Theorem 9 introduces new class of    generalizations, and it occurred unexpectedly that  \textit{Daykin--Eliezer--Carlitz result is not valid in necessary part for integrals!} Some more examples are  considered. The most general form of aithor's results are theorems 11 and 12.
In theorems 13, 14, 15 there are comparison results for different kinds of proved   generalizations.

In section 7 we consider applications of our results to Jackson $q$--integrals, Aszel inequality in Lorents spaces, discrete Fourier transforms, Minkowsky  inequality, two--sided estimates  for the Legendre complete integrals of the first kind, estimates for transmutations and integral operators, elasticity theory.

\newpage

 \UDC{517.162; 517.165}

 \begin{center}
 В книге:\\
 \bigskip
 Исследования по математическому анализу.\\
 \bigskip
 Серия: Математический Форум.\\
 \bigskip
 Редакторы:\\ Коробейник Юрий Федорович,\\
 Кусраев Анатолий Георгиевич.\\
 \bigskip
 Владикавказ:  Южный математический институт Владикавказского научного центра Российской Академии наук\\
 и Правительства Республики Северная Осетия-Алания ,\\
  2009.---T.~3.---С.~221--266.

 \bigskip
 \textbf{Уточнения и обобщения классических неравенств.}
 \bigskip

Автор: Ситник Сергей Михайлович,\\
доцент кафедры высшей математики\\
Воронежского института МВД, Воронеж, Россия.
\vskip 0.3cm

Электронный адрес: mathsms@yandex.ru
\vskip 0.3cm

Почтовый адрес: Ситник С.М.,\\
а.я. 12, Воронеж--5, Воронеж, 394005, Россия.

\vskip 0.3cm

\end{center}

 В  обзоре изложены уточнения и обобщения некоторых классических неравенств. Вначале рассмотрены уточнения неравенства Янга с их приложениями. Основное содержание обзора составляет изложение принадлежащего автору метода обобщений неравенств Коши--Буняковского с использованием средних значений. Рассмотрены обобщения для дискретного и интегрального случаев, их  варианты и сравнения между собой.  Приводятся приложения полученных результатов.

\bigskip
\epigraph{
Равенств не существует. Даже
в обычной человеческой жизни всегда встречается только неравенство.}{ Д. Митринович}

\section{1. Введение.}

Теория неравенств с приложениями отражена во многих  монографиях, см., например,~[3--4,6,9,15--16,23--24,30,32,35--36,60,62,72--73,75--79,87--88,92--93,
100--101,105,107--110]. Из них наряду с классическими трудами~[3,24,60,71,97--98,101,]  особенно отметим книги замечательного сербского математика Драгослава Митриновича ~[70--71,94--99,]. Его слова, взятые в качестве эпиграфа к настоящему обзору, отражают наблюдаемый уже несколько десятилетий устойчивый рост значения и роли неравенств в современной математике и её приложениях. Разнообразные
теоретические и численные оценки быстро вытесняют точные равенства, которые превращаются во многих случаях в прекрасные, но как бы случайные и вырожденные случаи при более полном и углублённом анализе решаемых задач.

В основаниях современной теории неравенств лежит сравнительно небольшое число классических результатов.
При этом каждое из основных знаменитых неравенств по мере своего существования обрастает многочисленными уточнениями, обобщениями и приложениями, которые способствуют более глубокому пониманию его сути и расширению набора возможных приложений.

Цель данного обзора --- изложение результатов автора по обобщениям неравенств Янга и Коши--Буняковского, в основном полученных в работах~[1,37--43] (см. также [8,51--57,106]), их приложениям к оценкам специальных функций [44,80--86] и ядер операторов преобразования [10--13,22,45--50]; при этом часть результатов публикуется впервые.

 \section{2. Уточнения неравенства Янга для двух чисел.}

Так сколько же неравенств на самом деле заключено в единственном на первый взгляд неравенстве Юнга (Янга)? И как это неравенство правильно называть? Неочевидные ответы на эти вопросы рассматриваются в этом пункте.

Прежде всего необходимо отметить, что по--русски правильно называть фамилию этого английского математика ---  Янг (William Henry Young). Традиция называть его на немецкий лад возникла из путаницы в переводе знаменитой монографии~[60], в которой вдобавок его фамилия в тексте дважды написана с ошибкой (Joung). Данное неравенство как и другое известное  неравенство Хаусдорфа--Янга для свёрток были доказаны Янгом с целью их применения в теории рядов Фурье. Аналогично принято неправильно произносить и фамилию его жены Грейс Чисхольм Янг, например, в названии известной  теоремы Янг--Данжуа--Сакса.

В 1912 году  В.~Янг доказал своё знаменитое неравенство в такой форме.

\Theorem { 1}
Пусть функция $f(x)$ является непрерывной и возрастающей на отрезке $~[0,c], c>0$. Пусть также $f(0)=0, a\in ~[0,c], b\in ~[0,f(c)]$. Тогда выполнено неравенство
$$
\int_0^af(x)\ dx + \int_0^bf^{-1}(x)\ dx \ge ab,\eqno(1)
$$
где $f^{-1}(x)$ есть обратная функция. При этом равенство в (1) достигается тогда и только тогда, когда $b=f(a)$.
\Endproc

Строго говоря, в первоначальной формулировке Янга функция была дифференцируемой, а не непрерывной, что было связано с методом его доказательства~[99].

Соотношение (1) геометрически выражает неравенство между площадями криволинейных трапеций, образованных графиками пары взаимно обратных функций. Аналитические доказательства, не использующие геометрические рассуждения, приведены, например, в~[3,60,99]. Современные доказательства как правило основаны на тождестве для левой части (1) при $b=f(a)$. Менее известно, что справедливо и обратное утверждение, доказанное японским математиком Т.~Такахаши в 1932 г.~[99].

Наиболее часто используется  неравенство Янга (1) при выборе $f(x)=x^{p-1}, p>1$. Тогда получаем в переменных $(x,y)$ числовое неравенство
$$
xy\le \frac{x^{p} }{p} +\frac{y^{q} }{q} \eqno(2)
$$
при условиях
$$
x\ge 0, y\ge 0,     p>1,   \frac{1}{p} +\frac{1}{q} =1. \eqno(3)
$$

Как и все классические неравенства, неравенство Янга может быть сформулировано с использованием средних значений. Тогда оно превращается в другое классическое неравенство: весовое среднее геометрическое не превосходит весового среднего арифметического:
 $$
 u^{a} v^{b} \le au+bv,\begin{array}{cc} {} & {} \end{array}a\ge 0,b\ge 0,a+b=1,u\ge 0,v\ge 0.
 $$

Принято считать, что круг вопросов, связанных с простейшим неравенством (2), исследован с исчерпывающей полнотой. Однако это не совсем так. Если обратить внимание на несимметричную форму правой части (2), то становится ясно, что {\bf {\it \underbar{неравенство Янга - это не одно, а пара неравенств}}}, причём второе из них, пропущенное  в (2), должно иметь такой вид:
$$
xy\le \frac{x^{q} }{q} +\frac{y^{p} }{p}. \eqno(4)
$$                                                                                             А вот тогда сразу возникает задача о сравнении неравенств (2) и (4), то есть об отыскании минимума правых частей из (2) и (4). Приведем численные примеры, которые показывают, что действительно между правыми частями (2) и (4) может существовать достаточно большой разброс.

\Primer \ 1: $x=5, y=130, p=4, q=4/3$;\\
$xy = 650,    \frac{x^{p} }{p} +\frac{y^{q} }{q} \approx 650,16502,$      $\frac{x^{q} }{q} +\frac{y^{p} }{p} \approx 71402508$.\\
В этом случае неравенство (2) лучше (на пять порядков!).

\Primer \ 2: $x=0,2, y=0,5, p=4, q=4/3$;\\
$xy = 0,1,   \frac{x^{p} }{p} +\frac{y^{q} }{q} \approx 0.29803,$     $\frac{x^{q} }{q} +\frac{y^{p} }{p} \approx 0,10334.$\\
А в этом случае неравенство (4) лучше (примерно в три раза).

Два приведенных примера иллюстрируют два типичных случая (хотя есть и третий!). Мы сразу сформулируем общий результат~[1,40]. Трудность была только в том, чтобы догадаться и найти данную оценку, после этого доказательство превращается в несложное упражнение.

Без ограничения общности далее будем предполагать, что выполнены условия  $y\ge x,$        $p\ge 2\ge q>1.$
\Theorem { 2}
Пусть выполнены условия  (3). Тогда

1. Если  $y\ge x\ge 1,$  то оценка (2) лучше, чем (4), то есть выполнены неравенства $$
xy\le \frac{[\min (x,y)]^{p} }{p} +\frac{[\max (x,y)]^{q} }{q} \le \frac{[\max (x,y)]^{p} }{p} +\frac{[\min (x,y)]^{q} }{q}. \eqno(5)
$$

2. Если  $1\ge y\ge x\ge 0,$  то оценка (4) лучше, чем (2), то есть выполнены неравенства
$$
xy\le \frac{[\max (x,y)]^{p} }{p} +\frac{[\min (x,y)]^{q} }{q} \le \frac{[\min (x,y)]^{p} }{p} +\frac{[\max (x,y)]^{q} }{q}. \eqno(6)
$$

3. Если   $y\ge 1\ge x\ge 0,$  то при данном $x$ существует единственное критическое значение $y=y_{cr}$, которое является решением трансцендентного уравнения
$$
\frac{x^{p} }{p} -\frac{x^{q} }{q} =\frac{y^{p} }{p} -\frac{y^{q} }{q}. \eqno(7)
$$
В этом случае при  $1\le y\le y_{cr}$  оценка (4) лучше, чем (2), то есть выполнены неравенства (6), а при  $y\ge y_{cr}$  оценка (2) лучше, чем (4), то есть выполнены неравенства (5).
\Endproc

Иначе говоря, если два числа лежат с одной стороны от единицы, то лучше одно из неравенств, с другой стороны---лучше другое, а если единица разделяет два числа, то реализуются оба случая.

Теперь приведем обещанный численный пример на третий случай в теореме 2.

\Primer \ 3: $x=0,5, p=4, q=4/3$. Тогда расчет, который мы опускаем (см.~[1,40]), дает критическое значение $y_{cr} \approx 1,35485.$
Выберем  $x=0,5, y=1,3<y_{cr}$; тогда получаем $xy =0,65$,
$$
\frac{x^{p} }{p} +\frac{y^{q} }{q} \approx 1,07973, \frac{x^{q} }{q} +\frac{y^{p} }{p} \approx 1,01166.
$$
Как и следует из теоремы 2, в этом случае лучше оценка (4).\\
А теперь пусть $x=0,5,\  y=1,4> y_{cr}$;  тогда получаем $xy=0,7$,
$$
\frac{x^{p} }{p} +\frac{y^{q} }{q} \approx 1.19025,\frac{x^{q} }{q} +\frac{y^{p} }{p} \approx 1.25804.
$$
Как и следует из теоремы 2, в этом случае лучше оценка (2).

Отметим, что наиболее полное изложение теории неравенства Янга и его обобщений приведено в~[99].

 \section{3. Дальнейшие обобщения и приложения неравенства Янга.}

Теперь наметим некоторые возможные  обобщения и приложения полученных результатов.

\subsec{1.} Как указал автору профессор В.~А.~Родин, можно рассмотреть более общий случай неравенств Янга (теперь мы знаем, что их два!) с парой произвольных взаимно дополнительных функций Янга. Так называют пару функций, удовлетворяющих соотношениям (обычно указывают только первое из них)
$$
xy\le N(x)+M(y) ,        xy\le N(y)+M(x) .                                               $$
Далее можно рассматривать разные случаи, от хороших функций до всё более плохих.

Случай 1 (самый хороший). Это исходное неравенство Янга (1), в котором
$$
N(x)=\int_0^xf(t)\ dt ,    M(y)=\int_0^yf^{-1}(s)\ ds.
$$
Тогда взаимно дополнительные функции Янга $N(x), M(y)$ являются и непрерывно дифференцируемыми, и выпуклыми, то есть хорошими.

Оказывается, что и в этой более общей формулировке можно полностью повторить результаты,  полученные выше  для простейшего случая. При этом роль пограничного условия  $x=1$  будет выполнять ненулевое решение некоторого трансцендентного уравнения с функцией $f(x)$.

Случай 2 (похуже). Пусть теперь дана функция $f(x), x\ge0$, неубывающая, неотрицательная, непрерывная справа при $x=0$ и удовлетворяющая условиям $f(0)=0, f(+\infty)=+\infty$. Введём функцию, называемую правой обратной для $f(x)$, по формуле
$$
g(s)=\sup_{f(x)\le s} x, s \ge 0.
$$
Тогда функции
$$
N(x)=\int_0^xf(t)\ dt ,    M(y)=\int_0^yg(s)\ ds
$$
образуют пару взаимно дополнительных функций Янга~[3,29,99]. Эти функции остались непрерывными и выпуклыми, но перестали быть непрерывно дифференцируемыми. Некоторый аналог теоремы 2 в этом случае ещё можно получить, так как осталась возможность брать производные.

Случай 3 (канонический). Кратко изложим необходимые для дальнейшего факты из выпуклого анализа~[17,21,33--34,58].

\Def{} Функция называется \textit{собственной}, если выполнено условие $f(x)>-\infty$, и хотя бы для одного $x$ из области определения $f(x)<+\infty$. \textit{Преобразованием Лежандра} собственной функции называется выражение
$$
(\mathcal{L}f(x))(y)=\sup_{x}(xy-f(x)).
$$
\Endproc

В многомерном случае произведение заменяется на скалярное произведение или другой линейный функционал. Для дифференцируемых функций можно выразить преобразование Лежандра неявно через градиент. Сама функция может быть фактически произвольной, но её преобразование Лежандра всегда будет выпуклой функцией с замкнутым надграфиком и полунепрерывной слева. Неравенства Янга в этой общей ситуации запишутся так:
$$
xy \le f(x)+(\mathcal{L}f(x))(y), xy \le f(y)+(\mathcal{L}f(y))(x),\eqno(8)
$$
и опять возникает задача об их сравнении. Но производной теперь нет, поэтому как их сравнивать---непонятно.

Для преобразования Лежандра произвольной собственной функции выполняется ещё одно красивое неравенство $\mathcal{L}(\mathcal{L}f)(x)\le f(x)$. Но можно заметить, что при подстановке в него конкретных функций, для которых преобразование Лежандра известно, всегда получается равенство. Это неспроста, и приводит к одному из самых важных фактов выпуклого анализа---теореме Фенхеля--Моро об инволютивности преобразования Лежандра. Она утверждает, что преобразование Лежандра инволютивно (то есть, его квадрат есть тождественное преобразование) тогда и только тогда, когда исходная  функция является собственной, выпуклой и её надграфик замкнут. Во всех рассмотренных выше случаях пары функций $N(x), M(x)$ удовлетворяли условиям этой теоремы, поэтому сами эти функции являлись преобразованиями Лежандра друг друга.
Возможно, неравенство (8) можно исследовать по нашей схеме, заменив применение производной на применение субдифференциалов~[18--19,58].

Случай 4 (самый плохой). В трёх приведённых случаях все функции были выпуклыми. Можно ли построить пару взаимно дополнительных по Янгу функций, хоть одна из которых  была бы невыпуклой? Оказывается, что да, можно. Этот случай рассматривается в теореме Бирнбаума--Орлича, приведённой в конце книги~[29]. В этой теореме приводится набор требований на одну из функций $N(x)$, которая может не быть выпуклой, чтобы вместе с её преобразованием Лежандра $M(x)$ они образовывали нужную пару и удовлетворяли неравенствам Янга. Как получить аналог теоремы 2 для этого случая совершенно неясно, так как теперь нет не только производных, но и пары инволютивных преобразований Лежандра.

\subsec{ 2.} Для большинства математиков неравенство Янга---это просто инструмент для доказательства неравенства Гёльдера. И это действительно одно из важнейших применений данного неравенства. Следует сказать прежде всего, что исторически справедливым было бы другое название: неравенство Роджерса--Гёльдера--Рисса.
Это неравенство было  впервые в эквивалентном виде доказано Леонардом Роджерсом в 1888 г., в другом эквивалентном виде Отто Гёльдером с указанием на вклад Роджерса  год спустя в 1889 г., а в привычном нам виде в 1910 г. Фридьешем (Фредериком?) Риссом в той же знаменитой работе, где им были введены пространства $l_p$ и $L_p$. К сожалению, работы замечательного английского математика Леонарда Джеймса Роджерса (1862--1933) были незаслуженно и надолго забыты. Приоритет Л.~Роджерса в доказательстве знаменитого неравенства был недавно восстановлен после публикаций польского математика Леха Малигранды~[88,91], хотя в хороших книгах это давно отмечалось, но почему-то в сносках~[3,60,99].

Итак, рассмотрим простейшее неравенство Роджерса--Гёльдера--Рисса для пары векторов с неотрицательными компонентами в конечномерном пространстве. Пусть
$$
a=(a_1,\dots a_n), b=(b_1,\dots b_n), A=\left(\sum_{k=1}^{n}a_{k}^p \right)^{\frac{1}{p}}, B=\left(\sum_{k=1}^{n}b_{k}^q \right)^{\frac{1}{q}},
$$
при условиях $p\ge 2\ge q>1$. Тогда обычное неравенство Р-Г-Р принимает вид
$$
\sum_{k=1}^{n}\frac{a_k}{A} \frac{b_k}{B} \le 1.\eqno(9)
$$
В стандартном доказательстве теперь нужно применить неравенство Янга к парам чисел $\frac{a_k}{A}, \frac{b_k}{B}$, которые очевидно все меньше единицы, что как раз и интересно, так как мы всегда окажемся в нестандартных условиях пункта 2 из теоремы 2. В результате получается

\Theorem { 3} Справедливо следующее уточнение дискретного неравенства Роджерса--Гёльдера--Рисса, следующее из теоремы 2:
$$
\sum_{k=1}^{n}\frac{a_k}{A} \frac{b_k}{B} \le
\sum_{k=1}^{n} \left(\frac{~[\max (\frac{a_k}{A},\frac{b_k}{B})]^{p} }{p} +\frac{~[\min (\frac{a_k}{A},\frac{b_k}{B})]^{q} }{q}\right) \le 1. \eqno(10)
$$
\Endproc

Аналогично может быть обобщено и соответствующее интегральное неравенство. Разумеется, любые уточнения неравенства Р-Г-Р в дискретном или интегральном случаях приводят к соответствующим уточнениям неравенства Минковского.

Отметим также разработанную А.~Г.~Кусраевым технику перенесения результатов для числовых неравенств
на равномерно полные числовые решётки и операторы в них.
Подобное распространение для неравенства Роджерса--Гёльдера--Рисса  проведено в~[89--90].
При этом основным инструментом становится представление  выпуклых положительных однородных функций в виде верхней огибающей семейства линейных функционалов, как в классических теоремах Минковского и Хёрмандера из выпуклого анализа.

\section{ 4. Классические и абстрактные средние значения.}

\subsec{ 4.1. Абстрактные средние.}
Средние величины первоначально возникли в Древней Греции при решении задач на пропорции в терминах отношения отрезков и других геометрических задач, см.~[9,109]. Так были введены простейшие средние: арифметическое, геометрическое, квадратичное и гармоническое
$$
A(x,y)=\frac{x+y}{2}, G(x,y)=\sqrt{x,y}, Q(x,y)=\sqrt{\frac{x^2+y^2}{2}}, x,y \ge 0,
$$
$$
H(x,y)=A(\frac{1}{x},\frac{1}{y})=\frac{2xy}{x+y}, \ x,y>0.
$$

Различные вопросы теории средних  достаточно подробно изложены в
литературе~[3,60,97,99]. Особенно отметим  специально
посвящённые этим вопросам монографии~[9,72,79,96,109], а также изложение в работах автора~[37--39].

Наблюдения за приведёнными выше простейшими средними приводят к выводу, что все они обладают некоторыми характерными общими свойствами: однородны, симметричны, монотонны, а также  при $x=y$ равны $x$. Это приводит к идее определять абстрактные средние при помощи соответствующих аксиом.

\sc   Определение. \rm
\textit{Абстрактным средним} двух неотрицательных чисел   называется
число $M(x,y)$, удовлетворяющее следующим аксиомам:

1) (свойство несмещённости)
$$
M(x,x) = x ,
$$

2) (свойство однородности)
$$
M(\lambda x, \lambda y) = \lambda \ M(x,y), \lambda > 0,
$$

3) (свойство монотонности по обоим аргументам)
$$
x_{2} > x_{1} \Rightarrow M(x_{2}, y) > M(x_{1}, y),\
 y_{2} > y_{1}\Rightarrow M(x, y_{2}) > M(x, y_{1});
$$

4)  (свойство симметричности)
$$
M(x,y) = M(y,x) .
$$

Очевидно, что если выполнено свойство симметричности, то тогда
достаточно монотонности лишь по одному аргументу. Из приведённых
выше  вытекает также важное свойство промежуточности
$$
\min(x,y)\le M(x,y) \le \max (x,y),
$$
которое часто включают в число основных (аналогично тому, как
включают в число аксиом нормы неотрицательность, хотя она следует из
других аксиом). Обычно подразумевается также непрерывность по обоим
аргументам, иногда для получения более детальных результатов нужны
дополнительно некоторая гладкость (например, если неравенства между
средними доказываются с использованием одной или нескольких
производных) и определённая модификация свойства аналитичности.

Для формулировки дальнейших результатов нам также потребуется
 понятие сопряженного среднего.

\sc Определение. \rm Сопряженным к абстрактному среднему
называется величина
$$
M^*(x,y) = \frac{xy}{M(x,y)};\  x,y>0. \eqno(11)
$$

Известны достаточно широкие классы средних, удовлетворяющих
приведенным аксиомам. Перечислим основные из них.

\subsec{ 4.2. Степенные средние.} Это самые известные средние, см., например,~[3,9,60,99]:

$$
M(x,y) = M_{\alpha} (x,y) = \left( \frac{x^{\alpha}+y^{\alpha}}{2}
\right)^{\frac{1}{\alpha}} , -\infty \le \alpha \le \infty \ ,
\alpha \neq 0\ ;\eqno(12)
$$
$$
M_{- \infty} (x,y) = \min (x, y)\ ,\  M_{0} = \sqrt{xy}\ , \
M_{\infty} (x, y) = \max (x, y)\ .
$$
Степенные средние образуют упорядоченную шкалу по параметру:
$$
\alpha_{1} > \alpha_{2} \Rightarrow M_{\alpha_1}(x,y) \ge
M_{\alpha_2} (x,y)\ ,\ \forall \ x, y\ .
$$
Три исключительных значения $\{\alpha=-\infty, 0, +\infty\} $ могут
быть получены из неисключительных предельным переходом. Аксиомы
абстрактного среднего 1)--4) проверяются непосредственно.

Кроме симметричных средних  используются и несимметричные. Самые
известные из них --- это  средние арифметическое и геометрическое с
неотрицательными весами $\alpha,\beta$ :
$$
A_{\alpha,\beta}(x,y)=\alpha x + \beta y,\\
G_{\alpha,\beta}(x,y)=x^{\alpha} y^{\beta},~\alpha +\beta=1. \eqno(13)
$$
Неравенство между ними---это  неравенство Янга, которое мы уже подробно рассмотрели.  Аналогично
вводится весовое среднее степенное произвольного порядка.

\subsec{ 4.3. Средние Т.~Радо.} Мы назовём так средние следующего вида:
$$
R_{\beta}(x,y) = \left( \frac{x^{\beta + 1} - y^{\beta + 1}}{(\beta
+ 1)(x-y)}\right) ^{\frac{1}{\beta}} , -\infty \le \beta \le
\infty ,\  \beta \neq 0, -1;\eqno(14)
$$
$$ R_{- \infty} (x,y) = \min (x, y),  R_{\infty} (x, y) = \max (x, y).
$$
Очевидно, что
$$
R_{-2} (x,y) = M_0 (x,y), R_{1} (x,y)=M_1 (x,y).
$$
Исключительные значения порождают пару гораздо менее известных средних:
логарифмическое
$$
R_{-1}(x,y)=L(x,y)=\frac{y-x}{\ln y-\ln x} \eqno(15)
$$
и то, которое в западной литературе почему--то называется identric, а автор предлагает называть многоэтажно--показательным или просто многоэтажным
$$
R_{0}(x,y)= \frac{1}{e} \left(
\frac{y^{y}}{x^{x}}\right)^{\frac{1}{y-x}}. \eqno(16)
$$

Средние Радо также образуют шкалу по параметру:
$$
\beta_{1} > \beta_{2} \Rightarrow R_{\beta_{1}} (x,y) \ge
R_{\beta_{2}} (x,y),~  \forall x, y;
$$
четыре исключительных значения  $\beta = \{-\infty, -1, 0,
+\infty\}$ могут быть получены из неисключительных предельным
переходом. Аксиомы абстрактного среднего  проверяются
непосредственно.

Конечно, величины $R_\beta$, которые являются интегральными средними от степенной функции, выписывались с момента возникновения Анализа. Но именно
Тибор Радо провел в 1935 г. в работе~[104] детальное исследование этих средних, которые понадобились ему для изучения субгармонических функций. В
частности‚ им были доказаны замечательные теоремы о связи двух
основных шкал средних  $M_\alpha$ и $R_\beta$. Так‚ из результатов
Радо следует‚ что средние из двух указанных шкал совпадают лишь в
следующих пяти случаях:
$$
M_{-\infty} = R_{-\infty}, M_{0} = R_{-2}, M_{\frac{1}{2}} =
R_{\frac{1}{2}}, M_{1} = R_{1}, M_{\infty} = R_{\infty}.
$$
Наиболее часто используемые средние ($\min$‚ $\max$‚ арифметическое‚
геометрическое) входят в обе шкалы. Другим замечательным результатом
Радо является полное описание множеств параметров  $(\alpha,\beta)$‚
при которых выполнены неравенства
$$
M_{\alpha_1}\le R_\beta \le M_{\alpha_2}, R_{\beta_1}\le M_\alpha
\le R_{\beta_2},\eqno(17)
$$

В частности‚ из его результатов получается‚ что фольклорное
неравенство для среднего логарифмического $M_{0} \le L \le M_{1}$,
которое  содержится еще в оригинальной работе В.\,Я. Буняковского
в качестве примера на приложения его интегрального неравенства,
может быть усилено до следующего:
$$
M_{0}(x,y) = \sqrt{xy}\le L(x,y) = \frac{x-y}{\ln x - \ln y}\le
$$
$$
\le M_{\frac{1}{3}}(x,y)=\left(
\frac{x^{\frac{1}{3}}+y^{\frac{1}{3}}}{2}\right)^{3},\eqno(18)
$$
причем порядки средних $0$ и $1/3$ являются неулучшаемыми. Последнее
неравенство многократно переоткрывалось в литературе.

К сожалению, условия  теорем Радо, при которых выполнены оценки (17),
приведены им не в явном виде. Так они цитируются и в известных монографиях, например,~[97,99]. Возможно, это одна из причин, по которой неравенства Радо многократно переоткрывались и переназывались. Расшифровка  условий оригинальной работы Радо не в
терминах некоторых дополнительных неравенств, а непосредственно
через параметры, получена в~[41]. Результаты содержит

\Theorem { 4} Справедливы следующие двусторонние неулучшаемые
оценки средних Радо через степенные средние:
$$
M_{\frac{\alpha + 2}{3}}\le R_{\alpha}\le M_{0}  \mbox{\ , при \ }
\alpha \in (-\infty, -2],
$$
$$
M_{0}\le R_{\alpha}\le M_{\frac{\alpha +2}{3}}  \mbox{\ , при \ }
\alpha \in ~[-2, -1],
$$
$$
M_{\frac{\alpha\ln2}{\ln(1+\alpha)}}\le R_{\alpha} \le
M_{\frac{\alpha +2}{3}} \mbox{\ , при \ }  \alpha \in (-1, -1/2],
$$
$$
M_{\frac{\alpha+2}{3}}\le R_{\alpha} \le M_{\frac{\alpha \ln
2}{\ln (1+\alpha)}} \mbox{\ , при \ } \alpha \in ~[-1/2, 1),\\
$$
\mbox{(при $\alpha = 0$  последнее неравенство понимается в
предельном смысле}\\ $M_{\frac{2}{3}}\le R_0 \le M_{\ln 2}$),
$$M_{\frac{\alpha\ln
2}{\ln (1+\alpha)}}\le R_\alpha \le M_{\frac{\alpha +2}{3}}
\mbox{\ , при \ } \alpha \in ~[1,\infty].\\
$$
\Endproc

Данная теорема была  распространена на случай более общих средних
Джини--Радо (Столярского) в работах~[102--103].

Отметим также полученное в~[42]
неравенство для среднего логарифмического
$$
L\left( M_{\frac{1}{2}}(x,y), M_0(x,y)\right) \le L(x,y) \le
L(M_1((x,y), M_0(x,y)),
$$
в котором порядки средних наилучшие.

Неравенства между средними --- это различные способы их
упорядочивания. Иногда используются достаточно экзотические варианты
такого упорядочивания, например, с использованием преобразования
Фурье и положительно--определенных функций, которые тем не
менее находят полезные приложения, см.~[79].

\subsec{ 4.4. Средние Джини и Лемера.}

Итальянский статистик Коррадо Джини ввёл названные его именем средние с двумя параметрами в 1938 г.  по формулам
$$
Gi_{u,v}(x,y)=\left( \frac{x^u+y^u}{x^v+y^v}\right)^{\frac{1}{u-v}}, u\neq v,
\eqno(19)
$$
$$
Gi_{u,v}(x,y)=\exp \left(\frac{x^u\ln x+y^u \ln y}{x^u+y^u}\right), u=v\neq0,
$$
$$
Gi_{u,v}(x,y)=G(x,y), u=v=0.
$$

Основным предметом изучения Джини было имущественное неравенство, для его описания он и ввёл индекс Джини, основанный на предыдущих формулах. Это одно из основных понятий современной социальной статистики.

Важный частный случай средних Джини получается из (19) при $v=u-1$.
Эти средние были переоткрыты в  1971 г. Д.~Лемером и имеют вид
$$
Le_u(x,y)=\frac{x^{u+1}+y^{u+1}}{x^u+y^u}. \eqno(20)
$$

\subsec{ 4.5. Квазиарифметические средние.}

\sc Определение. \rm Пусть дана неотрицательная монотонная функция $f(x)$, набор неотрицательных чисел $x=(x_1,x_2,\cdots,x_n)$ и неотрицательных весов
$p=(p_1,p_2,\cdots,p_n)$. \textit{Квазиарифметическим средним} называется выражение
$$
Q_p(x)=f^{-1}\left( \sum_{k=1}^n p_k f(x_k)\right), \sum_{k=1}^n p_k=1.
$$
В частном случае $f(x)=x$ получаем обычное весовое среднее арифметическое, чем объясняется название. Это среднее стало широко известно после выхода монографии~[60], в которой оно применялось для изучения обобщённой выпуклости относительно пары функций. Основные результаты для квазиарифметического среднего были получены в работах Колмогорова, Нагумо и Де Финетти. Они в существенном сводятся к доказательству двух основных свойств:

1. Квазиарифметические средние совпадают тогда и только тогда, когда порождающие их функции связаны линейным соотношением.

2. Непрерывное  квазиарифметическое среднее однородно тогда и только тогда, когда оно совпадает со средним степенным.

Данные средние встречаются в огромном числе задач.  Например, они используются в недавних работах акад. В.~П.~Маслова по применению нелинейных моделей в математической экономике~[25--28]. В монографии~[61] квазиарифметические  средние применяются для детального анализа различных видов энтропии.

\subsec{ 4.6. Итерационные средние.}

Рассмотрим итерационный процесс при
заданных стартовых значениях $x_0, y_0$  и паре
абстрактных средних $(M, N)$\ : \\
$$
x_{n+1} = M(x_n, y_n), \ y_{n+1} = N(x_n, y_n).
$$

В общем случае получаем некоторую динамическую систему на плоскости,
ко\-то\-рая имеет интересное асим\-пто\-ти\-чес\-кое по\-ве\-де\-ние
(динамику) при $n\rightarrow~\infty$, но очень сло\-жна для изучения
даже при простейшем выборе пары средних $(M, N)$.

\sc Определение. \rm Пусть существует общий предел
последовательностей $x_n$ и $y_n$. Тогда он называется
\textit{итерационным средним} и обозначается
$$
\mu (M,N|\ x_0, y_0) = \mu(x_0, y_0) = \lim x_n = \lim y_n. \eqno(21)
$$

Непосредственно проверяется, что итерационное среднее
пары абстрактных средних также является новым абстрактным средним и
наследует соответствующие аксиомы 1) - 4). Также принято другое обозначение вместо $\mu(M,N)$, а именно $M\otimes N$.

Самое известным итерационным средним является арифметико--геометрическое (AGM---arithmetic--geometric mean), изученное Лежандром и Гауссом. Оно получается при выборе $M=M_1$, $N=M_0$ \   и выражается по формуле
$$
\mu(M,N|\ x_0,y_0) = \frac{\pi}{2} \frac{x_0}{K\left(\sqrt{1-\left(
\frac{y_0}{x_0}\right)^{2}}\right)}, 0 < y_0 <x_0, \eqno(22)
$$
где $K(x)$ --- полный эллиптический интеграл Лежандра первого рода.

Разумеется, можно рассматривать и огромное число итераций с другими средними, или менять их по ходу итераций, по аналогии с методом Зейделя.

Теория итерационных средних в настоящее время активно развивается и
является, в частности, источником возникновения новых классов
специальных функций --- гипертрансцендентных функций. Значительная часть этой теории изложена в книгах Джонатана и Питера Борвейнов с соавторами~[65--69,5].

В заключение перечисления различных средних отметим, что некоторые
элементарные операции не выводят из этого класса. Например, среднее от
пары средних---это вновь среднее. Данную процедуру можно неограниченно продолжать.  Несмотря на свою очевидность, это
формулировка очень мощного метода, порождающего в конкретных случаях
практически все известные процедуры построения новых средних из
простейших.

Полезным необычным средним является также понятие медианты двух
дробей, вводимое в теории рядов Фарея, по которому неграмотные школьники и студенты и складывают дроби. Далее мы используем это
среднее.

Таким образом, можно сделать основной для данного раздела вывод, что
\textit{существует значительное число конкретных примеров
абстрактных средних}. Они используются далее в качестве
"кирпичиков"\,, из которых, например, собираются различные обобщения неравенств
Коши--Буняковского в интегральном и дискретном случаях.

Необходимо подчеркнуть фундаментальную роль средних степенных и
Радо, так как все остальные перечисленные выше известные средние
выражаются через них. Поэтому именно эти две шкалы средних играют
основную роль в данной работе.

\subsec{ 4.7. Полное описание функций средних и энтропии для двух переменных.}

Итак, средних очень много. Тем удивительнее, что удалось получить описание всех "хороших"\ средних от двух переменных. В этом пункте излагаются соответствующие результаты, полученные автором.

\Theorem { 5} Произвольное абстрактное среднее, обладающее свойствами несмещённости, однородности, симметричности,
монотонности и непрерывности по обоим аргументам представляется в виде
$$
M (x,y)= (x+y)h(\ln \frac{y}{x}),  \eqno(23)
$$
 где  $h(t)$ является определенной на всей оси непрерывной
 четной функцией, удовлетворяющей следующим условиям при $t\ge
0$:
$$
h(0)=\frac{1}{2} ,\\
$$
$$
\frac{e^{t_1}(e^{t_2}+1)}{e^{t_2}(e^{t_1}+1)}\ \le \
\frac{h(t_1)}{h(t_2)} \le  \frac{(e^{t_2}+1)}{(e^{t_1}+1)}, \
t_1\le t_2.
$$
Справедливо и обратное: каждой функции $h(t)$ с указанным
набором свойств соответствует по формуле (23) несмещённое,
однородное, симметричное, монотонное и непрерывное по обоим
аргументам среднее.
\Endproc

Введенная функция $h(t)$ связана с функцией однородности
соответствующего среднего $M (x,y)$.

Полученное полное описание средних значений позволяет дать
исчерпывающее описание функций энтропии от двух переменных, что
является важным для теории информации и других приложений.

\sc Определение. \rm Обобщённой энтропией, соответствующей произвольному абстрактному среднему, называется величина
$$
E(x,y)=-\ln( M( x,y )).  \eqno(24)
$$

Это определение не надумано, а отражает известные конкретные функции энтропии, которые общеприняты в термодинамике и информатике.

Справедлива и обратная формула, выражающая среднее через произвольную энтропию
$$
M( x,y)=\exp(-E(x,y)).
$$

Выбор конкретных средних в определении (24) порождает все известные функции энтропии. Так при выборе весового среднего геометрического получается определение энтропии по Шеннону, а при выборе  среднего Джини получается определение энтропии по Реньи~[61,108]:
$$
E_1(x,y)=-p_1\ln x - p_2\ln y,
E_2(x,y)=-\frac{1}{\alpha -1}\ln \frac{(p_1 x^\alpha + p_2 y^\alpha)}{p_1 x+p_2 y}.
$$

Из предыдущей теоремы сразу получаем полное описание всех функций энтропии от двух переменных.

\Theorem { 6} Произвольная функция энтропии определяется по формуле
$$
E(x,y)=-\ln(x+y)-\ln h(\ln \frac{y}{x}),
$$
с некоторой функцией  $h(t)$, где  $h(t)$ является определенной на всей оси непрерывной  четной функцией, удовлетворяющей следующим условиям при $t\ge
0$:
$$
h(0)=\frac{1}{2} ,\\
$$
$$
\frac{e^{t_1}(e^{t_2}+1)}{e^{t_2}(e^{t_1}+1)}\ \le \
\frac{h(t_1)}{h(t_2)} \le  \frac{(e^{t_2}+1)}{(e^{t_1}+1)}, \
t_1\le t_2.
$$
\Endproc

\section{5. Обобщения дискретных неравенств Коши--Буняковского.}

\subsec{ 5.1. Исторические сведения.}

Неравенство Коши--Буняковского является классическим неравенством
Анализа. Для конечных сумм и рядов оно было доказано О.~Коши  в
1821 году, а для интегралов --- В.~Я.~Буняковским в 1859 году.
Отметим, что в 2004 году в Киеве состоялась международная
конференция, посвященная двухсотлетию со дня рождения академика
Виктора Яковлевича Буняковского, на которой автор выступал  с пленарным
докладом по тематике данной статьи.

Курс Анализа Коши, содержащий данное неравенство, вышел в России в 1831 г. в переводе Буняковского.  Собственная работа Буняковского содержит много примеров на применение его неравенства и вычислений, но доказательства там нет. Всем понятно, что и Коши, и Буняковский его знали.

На самом деле, неравенство Коши сразу следует из тождества Лагранжа, которое было доказано им почти за 50 лет до Коши при $n=3$, а при $n=2$ было известно ещё Диофанту. Лагранж вывел своё тождество при изучении геометрии пирамид.

Интегральное неравенство Буняковского было переоткрыто через
двадцать пять лет Г.~А.~Шварцем в 1884. Принято считать, что им впервые дана теперь привычная всем формулировка  в терминах скалярного
произведения, что должно оправдать его долю в наименовании этого неравенства, но это неправда. На самом деле большая работа Шварца была посвящена его альтернирующему методу. При оценке сходимости  он использовал соответствующее неравенство, которому в  работе посвящён один абзац из 6 строк. В нём дан классический вывод через квадратный трёхчлен.  Кстати, обдумывание этой работы вскоре привело Эмиля Пикара к первой версии метода последовательных приближений для дифференциальных уравнений. Герман Амандус Шварц многое открыл сам, и нет никакой нужды приписывать ему чужие достижения. Вспомним альтернирующий метод Шварца, теорему о равенстве смешанных производных, решение задачи Фаньяно, интеграл Кристоффеля--Шварца, принцип симметрии Римана--Шварца, сапог Шварца, производную Шварца и т.д.

Идея обобщения неравенства на пространства со скалярным произведением была предвосхищена Г.~Грассманом, а в привычном нам виде результат был опубликован Г.~Вейлем в 1918 г. Первым включил неравенство Коши--Буняковского в свои монографии Фон Нейман.

Сейчас издано несколько монографий, целиком посвящённых рассмотрению неравенства Коши--Буняковского~[75--77,107,110], особо из них отметим недавнюю замечательную книгу Стила и первую из монографий Севера Драгомира. Там есть достаточно полная библиография, поэтому мы не будем цитировать
многочисленные работы, посвящённые различным обобщениям данного
неравенства.

\subsec{ 5.2. Неравенства Милна и Колбо.}

Итак, рассмотрим дискретное неравенство Коши--Буняковского
$$
\left(\sum_{k=1}^nx_k\cdot y_k\right)^2 \le \left(\sum_{k=1}^nx_k^2\right)\cdot\left(\sum_{k=1}^ny_k^2\right).\eqno(25)
$$

В [60,75,99] приведено его интересное обобщение  в форме
$$
\left(\sum_{k=1}^nx_k\cdot y_k\right)^2 \le
\left(\sum_{k=1}^n(x_k^2+y_k^2)\right)\cdot\left(
\sum_{k=1}^n\frac{x_k^2 y_k^2}{x_k^2+y_k^2}\right)
\le\eqno(26)
$$
$$
\le\left(\sum_{k=1}^nx_k^2\right)\cdot\left(\sum_{k=1}^ny_k^2\right).
$$

Это  известное  неравенство,  которое было доказано Милном
в 1925 году при исследовании задачи из астрономии
по вычислению коэффициента звездного поглощения, и с тех пор
неоднократно переоткрывалось. Неравенство
Милна явилось  одной из отправных точек для получения изложенных далее результатов по обобщению неравенства Коши--Буняковского. В~[60,75] приведены также некоторые подобные оценки. Интересные приложения неравенства Милна к теории электрических цепей и некоторые обобщения приведены в~[105]. Установлена связь этого неравенства с физическим принципом Максвелла.  Только неравенство в тексте книги названо чужим именем, но в заключительном замечании Сандор восстановил справедливость, заодно перечислив ещё ряд переоткрытий этого результата.

Другим подобным результатом является неравенство, доказанное Колбо в  1965 г.
$$
\left(\sum_{k=1}^nx_k\cdot y_k\right)^2 \le
\sum_{k=1}^nx_k^{1+\alpha} y_k^{1-\alpha}
\sum_{k=1}^nx_k^{1-\alpha} y_k^{1+\alpha}\le
\left(\sum_{k=1}^nx_k^2\right)\cdot\left(\sum_{k=1}^ny_k^2\right),
$$
справедливое при $\alpha\in[0,1]$.

\subsec{ 5.3. Теорема Дэйкина--Элиезера--Карлица.}

Достаточно очевидно, что приведённые выше обобщения неравенства Коши--Буняковского, принадлежащие Милну и Колбо, имеют общую природу. Она была выявлена в следующем результате~[75,99].

\Theorem { Дэйкина--Элиезера--Карлица (\textbf{CDE})}
Для того, чтобы для любых двух последовательностей неотрицательных чисел было выполнено неравенство
$$
\left(\sum_{k=1}^nx_k\cdot y_k\right)^2 \le
\sum_{k=1}^n f(x_k,y_k)\sum_{k=1}^n g(x_k,y_k)\le\eqno(27)
$$
$$
\left(\sum_{k=1}^nx_k^2\right)\cdot\left(\sum_{k=1}^ny_k^2\right),
$$
необходимо и достаточно, чтобы пара функций $f(x,y),g(x,y)$ удовлетворяла следующим трём условиям
$$
1) \  f(x,y) g(x,y)=x^2 y^2,\eqno(28)
$$
$$
2) \ f(\lambda x,\lambda y)=\lambda^2 f(x,y),
$$
$$
3) \ \frac{y f(x,1)}{x f(y,1)}+\frac{x f(y,1)}{y f(x,1)}\le \frac{x}{y}+\frac{y}{x}. \eqno(29)
$$
\Endproc

Неравенство Милна является частным случаем теоремы CDE при выборе $f(x,y)=x^2+y^2$, а неравенство Колбо при выборе $f(x,y)=x^{1+\alpha} y^{1-\alpha}, \alpha\in[0,1]$.

Если проанализировать доказательство теоремы CDE, то становится ясным, что в достаточной части авторы переоткрывают аналог тождества Лагранжа, который существует для неравенства Чебышёва. Поэтому доказательство достаточности лучше сразу провести с применением неравенства Чебышёва, что и проще и не затемняет сути дела. Замечательным является найденный авторами способ доказательства необходимости теоремы и обнаруженное ими условие (29). Это гибридное условие на функцию, записываемое единой простой формулой, из которого следует одновременно её однородность и монотонность, аналогично тому как из известного условия Годуновой--Левина следует одновременно монотонность и выпуклость. На мой взгляд роль таких гибридных условий, объединяющих в одной формуле несколько стандартных, пока в Анализе недооценена, их значение ещё предстоит раскрыть.

Теперь сделаем следующий важный шаг. Перепишем неравенства Милна и Колбо в терминах средних значений. Начнём с первого из них.
$$
\left(\sum_{k=1}^nx_k\cdot y_k\right)^2 \le
\left(\sum_{k=1}^n {(M_2(x_k,y_k))}^2\right)\cdot\
\left(\sum_{k=1}^n {({M_2^*}(x_k,y_k))}^2\right)
\le\eqno(30)
$$
$$
\le\left(\sum_{k=1}^nx_k^2\right)\cdot\left(\sum_{k=1}^ny_k^2\right),
$$
где $M_2=Q$---среднее квадратичное, а $M_2^*=M_{-2}$---его сопряжённое (11). Это предлагаемая мной запись неравенства Милна через среднее квадратичное.

Далее, запишем
$$
\left(\sum_{k=1}^nx_k\cdot y_k\right)^2 \le
\left(\sum_{k=1}^n {(G_{\frac{1+\alpha}{2},\frac{1-\alpha}{2}}(x_k,y_k))}^2\right)\cdot\
$$
$$
\cdot\left(\sum_{k=1}^n {({G^*}_{\frac{1+\alpha}{2},\frac{1-\alpha}{2}}(x_k,y_k))}^2\right)
\le\left(\sum_{k=1}^nx_k^2\right)\cdot\left(\sum_{k=1}^ny_k^2\right),\eqno(31)
$$
где $\alpha\in[0,1]$, $G_{\frac{1+\alpha}{2},\frac{1-\alpha}{2}}$---весовое среднее геометрическое (13), ${G^*}_{\frac{1+\alpha}{2},\frac{1-\alpha}{2}}=
G_{\frac{1-\alpha}{2},\frac{1+\alpha}{2}}$---его сопряжённое (11). Это предлагаемая мной запись неравенства Колбо через весовое среднее геометрическое.

Теперь мы готовы записать теорему Дэйкина--Элиезера--Карлица в терминах средних, начав соединение основных тем данного обзора.

\Theorem { 7 (CDE в терминах средних значений)}
Пусть дано абстрактное несмещённое, однородное и монотонное по обоим аргументам среднее $M(x,y)$, $M^*(x,y)$---его сопряжённое среднее (11). Тогда справедливо следующее обобщение неравенства Коши--Буняковского
$$
\left(\sum_{k=1}^nx_k\cdot y_k\right)^2 \le
\left(\sum_{k=1}^n {(M(x_k,y_k))}^2\right)\cdot\
\left(\sum_{k=1}^n {({M^*}(x_k,y_k))}^2\right)
\le\eqno(32)
$$
$$
\le\left(\sum_{k=1}^nx_k^2\right)\cdot\left(\sum_{k=1}^ny_k^2\right).
$$
Верно и обратное: любое обобщение неравенства Коши--Буняковского из теоремы
Дэйкина--Элиезера--Карлица (CDE) в форме (27) может быть записано в виде (32) с некоторым средним, удовлетворяющим перечисленным выше свойствам.
\Endproc

Отметим, что условие симметричности среднего в этой теореме не требуется, что допускает такие случаи, как неравенство Колбо.

По моему мнению, переформулировка теоремы
CDE в терминах средних делает ее гораздо более понятной и
прозрачной, а также снабжает эту теорему колоссальным числом конкретных
примеров из теории средних величин в дополнение к скромному набору
из двух примеров (неравенств Милна и Колбо) во всех  известных монографиях.
Становится понятным, зачем мы описывали огромное множество средних в части 4 данного обзора --- с их помощью теперь можно повыписывать также огромное число обобщений неравенства Коши--Буняковского совершенно различного вида. Важной задачей является также сравнение этих обобщений между собой, когда это возможно.

\section{6. Обобщения интегральных неравенств Коши--Буняковского.}

\subsec{ 6.1. Обобщения с различными средними.}

Мы будем рассматривать обобщения  интегрального неравенства Коши-Буняковского
следующего вида:
$$
  \left(\int^b_a f(x)g(x)\,dx\right)^2 \le
\int^b_a\Phi_1(f,g)\,dx\cdot
 \int^b_a\Phi_2(f,g)\,dx \le\eqno(33)
$$
$$
 \le\int^b_a(f(x))^2\,dx \cdot
\int^b_a(g(x))^2\,dx , \phantom{aaaa}
$$
которые должны выполняться для достаточно
произвольных функций $f(x)$, $g(x)$ при выборе некоторых
неотрицательных функций $\Phi_1$, $\Phi_2$.

Будем считать далее выполненными следующие простейшие предположения:
все встречающиеся функции непрерывны, интегралы существуют в смысле
Римана, промежуток интегрирования конечен, переменная одна. От
каждого из этих ограничений можно отказаться и рассмотреть
соответствующие обобщения (интегрируемые в определенном смысле
функции, интеграл Лебега, неограниченное множество интегрирования,
функции  нескольких переменных), но эти вопросы здесь не
рассматриваются.

\Theorem { 8}
Пусть  $M$ - произвольное несмещённое,
однородное, монотонное (необязательно симметричное!) абстрактное
среднее, $M^*$ - сопряжённое к нему. Тогда справедливо обобщение
интегрального неравенства Коши - Буняковского вида (33) при выборе
$$
\Phi_1(f,g)=(M(f,g))^2, \Phi_2(f,g)=(M^*(f,g))^2.\eqno(34)
$$
\Endproc

Таким образом, мы получаем сформулированный в терминах средних
результат, что справедлив интегральный аналог \textit{достаточной}
части  теоремы CDE Карлица--Дэйкина--Элиезера  для дискретного случая.

Сделаем важное для дальнейшего замечание, что левое неравенство из (33) при нашем выборе (34) в силу определения сопряжённого среднего (11) очевидно,
так как оно само сводится к обычному неравенству Коши--Буняковского. Таким образом,
содержательным в теореме 8 является правое неравенство между
произведениями интегралов. Однако, и левое очевидное неравенство из
(33) может быть с успехом использовано для самоуточнения неравенства Коши--Буняковского, см. далее.

 В частном случае при выборе степенного среднего
 $M(f, g)=M_{\alpha}(f,g)$  получаем

\sc Следствие 8.1. \rm Для неотрицательных непрерывных функций
$f(x)$, $g(x)$, одновременно не обращающихся в ноль на $[a, b]$,
справедливо следующее обобщение неравенства Коши--Буняковского:
$$
\left(\int^b_af(x)g(x)dx\right)^2\le
\int^b_a[M_\alpha(f,g)]^2dx\cdot\int^b_a[M_{-\alpha}(f,g)]^2dx=\eqno(35)
$$
$$
=\int^b_a(f^\alpha+g^\alpha)^{2/\alpha}dx\cdot\int^b_a
f^2g^2(f^\alpha+g^\alpha)^{-2/\alpha}dx\le\int^b_af^2dx\cdot\int^b_ag^2dx.
$$

Выпишем наиболее эффектные из полученных неравенств с частными
случаями степенных средних  явно. Они соответствуют выбору среднего
полуквадратичного, арифметического и квадратичного:

\begin{eqnarray*}
\left(\int^b_af(x)g(x)\,dx\right)^2 \le \int^b_a
{\left(\sqrt{f(x)}+\sqrt{g(x)}\right)}^4\,dx \ \cdot &\\
\cdot \int^b_af^2g^2 /{\left(\sqrt{f(x)}+\sqrt{g(x)}\right)}^4\,dx
\le  \int_a^b f^2(x)\,dx \cdot \int_a^b g^2(x)\,dx,&
&\\
\left(\int^b_af(x)g(x)\,dx\right)^2 \le \int^b_a
{\left(f(x)+g(x)\right)}^2\,dx \ \cdot &\\
\cdot \int^b_af^2g^2 /{{\left(f(x)+g(x)\right)}^2}\,dx \le
\int_a^b f^2(x)\,dx \cdot \int_a^b g^2(x)\,dx,&
&\\
\left(\int^b_af(x)g(x)\,dx\right)^2 \le \int^b_a
\left(f^2(x)+g^2(x)\right)\,dx \ \cdot &\\
\cdot \int^b_af^2g^2 /\left(f^2(x)+g^2(x)\right) ,dx \le \int_a^b
f^2(x)\,dx \cdot \int_a^b g^2(x)\,dx.& &
\end{eqnarray*}
При $\alpha=2$ получаем  интегральный вариант неравенства Милна.
А при  $\alpha = + \infty $ получается такое занятное неравенство.

\sc Следствие 8.2. \rm Справедливо следующее уточнение
неравенства Коши--Буняковского
$$
\left(\int^b_af(x)g(x)dx\right)^2\le\int^b_a[\max(f,g)]^2dx\cdot
\int^b_a[\min(f,g)]^2dx\le\eqno(36)
$$
$$
\le\int^b_af^2(x)dx\cdot \int^b_ag^2(x)dx.
\phantom{aaaaaaaaaaa}
$$

При первом взгляде на это неравенство кажется, что средняя часть
должна сводится к одной из крайних. Однако, это не так. То, что все
три части в (36) могут быть различны, доказывает пример "конвертика"
\ : $a=0$, $b=1$, $f(x)=x$, $g(x)=1-x$.
Вычисления показывают, что в этом случае (36) сводится к
неравенствам $\frac{1}{36}< \frac{7}{144}< \frac{1}{9}$. \smallskip
Кстати, ссылки, где это
простое неравенство написано раньше, автор не знает.

 Проанализируем набор полученных неравенств со средними степенными (35) в
зависимости от параметра $\alpha$. Прежде всего, в силу симметрии в
(35) можно ограничиться значениями $\alpha\ge 0$. При  $\alpha = 0$
средняя часть (35), которая понимается в предельном смысле, равна
левой. Поэтому, при $\alpha \approx 0$ будем получать хорошие
приближения к левой части. Отметим интересную особенность: при
$\alpha\rightarrow\infty$ средняя часть в (35) стремится к средней
части в рассмотренном выше неравенстве (35)‚ а не к правой части
(35). Таким образом‚ в (35) возникает "зазор"\  с правой частью,
который‚ возможно‚ может быть заполнен оценками другого вида.

При выборе $M(f, g)=R_{\beta}(f,g)$ получаем
\smallskip

\sc Следствие 8.3. \rm Справедливо следующее уточнение
неравенства Коши-Буняковского в терминах средних Радо :
$$
\left(\int^b_af(x)g(x)dx\right)^2\le\int^b_a\left[R_\beta(f,g)\right]^2dx
\cdot\int^b_a\left[\frac{1}{R_\beta\left(\frac{1}{f},
\frac{1}{g}\right)}\right]^2dx \le \eqno(37)
$$
$$
\le\int^b_af^2(x)dx\int^b_ag^2(x)dx.
\phantom{aaaaaaaaaaaaaaaa}
$$

Выпишем наиболее эффектные из полученных неравенств (37) с частными
случаями средних Радо  явно. Они включают  выбор среднего
логарифмического и многоэтажно--показательного.

\begin{eqnarray*}
\left(\int^b_af(x)g(x)dx \right)^2 & \le
\int^b_a\left(\frac{f-g}{\ln\frac{f}{g}}\right)^2\,dx \cdot
\int^b_af^2g^2/\left(\frac{f-g}{\ln\frac{f}{g}}\right)^2\,dx& \le  \\
&\le  \int_a^b f^2(x)\,dx \cdot \int_a^b g^2(x)\,dx,&
\\
\left(\int^b_af(x)g(x)dx \right)^2 & \le
\int^b_a\left[\frac{f^f}{g^g} \right]^{\frac{2}{f-g}}\,dx \cdot
\int^b_af^2g^2/\left(\frac{f^f}{g^g}
\right)^{\frac{2}{f-g}}\,dx& \le  \\
&\le \int_a^b f^2(x)\,dx \cdot \int_a^b g^2(x)\,dx,&
\\
\left(\int^b_af(x)g(x)dx\right)^2& \le
\int^b_a\left(f^2+fg+g^2\right)\,dx
\cdot\int^b_a\frac{f^2g^2}{f^2+fg+g^2}\,dx& \le \\
&\le \int_a^b f^2(x)\,dx \cdot \int_a^b g^2(x)\,dx.&
\end{eqnarray*}

Анализ полученных неравенств со средними  Радо  в зависимости от
параметра $\beta$ аналогичен проведённому выше для степенных
средних. Отметим и здесь наличие "зазора" \   с правой частью.
\smallskip

При выборе $M(f,g)=\mu(M_1, M_0|f,g)$ получаем
\smallskip

\sc Следствие 8.4. \rm Справедливо следующее уточнение неравенства
Коши - Буняковского:
 \begin{eqnarray*}
 \left(
\int^b_af(x)g(x)dx\right)^2\le\int^b_a{\left[\frac{\max(f,g)}
{K\left(\sqrt{1-\left(\frac{\min(f,g)}
{\max(f,g)}\right)^2}\right)}\right]}^2dx \cdot\\
\cdot\int^b_a ( \min(f,g) )^2 \Biggl(
K\left(\sqrt{1-\left(\frac{\min(f,g)} {\max(f,g)}\right)^2}\right)
\Biggr) ^2\,dx \le\int^b_af^2\,dx\int^b_ag^2\,dx,
\end{eqnarray*}
где $K(x)$ есть полный эллиптический интеграл Лежандра 1 рода.

Отметим совершенно экзотический характер последнего неравенства: это
неравенство между \textit{произвольными} функциями, но которые стоят
под знаком \textit{конкретной специальной} функции - эллиптического
интеграла Лежандра! Вместе с тем неравенства, в которых произвольные
функции стоят под знаками \textit{элементарных} функций, например,
экспоненты или логарифма, широко используются в математике. Так,
подобный вид имеют неравенства Джона--Ниренберга в теории
пространств ВМО и неравенства Зигмунда в гармоническом анализе.

Чтобы  продемонстрировать возможность применения несимметричных средних, приведем соответствующие результаты.

\sc Следствие 8.5.\rm Для неотрицательных непрерывных функций
$f(x)$, $g(x)$ одновременно не обращающихся в
ноль на $[a, b]$, при $0\le\alpha \le 1$   справедливы
следующие   обобщения  неравенства  Коши--Буняковского:
\begin{eqnarray*}
\left(\int^b_af(x)g(x)\,dx \right)^2\le
\int^b_a[\alpha f(x)+(1-\alpha)\ g(x)]^2\,dx \ \cdot &\nonumber\\
\cdot\int^b_af^2(x)g^2(x)/[\alpha f(x) + (1-\alpha)\ g(x)]^2\,dx
\le
\int^b_af^2(x)\,dx \int^b_a g^2(x)\,dx , &\\
\left(\int^b_af(x)g(x)\,dx \right)^2\le
\int^b_a[f^\alpha(x)g^{(1-\alpha)}(x)]^{2}\,dx
\int^b_a [f^{(1-\alpha)}(x)g^\alpha(x)]^{2}\,dx \le \nonumber\\
\le \int^b_af^2\,dx  \int^b_ag^2\,dx.
\phantom{aaaaaaaaaaaaaaaaaaaaaaaaaa}
\end{eqnarray*}

Последнее неравенство  есть интегральный аналог результата
Колбо.

\subsec{ 6.2. Обобщения с логарифмическими производными.}

 Для дифференцируемых функций получен другой
ряд обобщений. Они получают особенно важное значение при сравнении
наших интегральных неравенств и теоремы CDE для дискретного случая.

\Theorem { 9}
Для подходящих функций справедливо обобщение
интегрального неравенства Коши - Буняковского вида (33) при
следующем выборе:
$$
 \Phi_1(f,g)=\exp\left(2\int^x_aM(Lf(y), Lg(y))\,dy\right),\
\Phi_2(f,g)=\frac{f^2g^2}{\Phi_1(f,g)} ,\eqno(38)
$$
где ~$Lh(t)$ - логарифмическая производная $Lh=\frac{h'}{h}$, $M$ -
произвольное несмещённое, однородное, монотонное (необязательно
симметричное!) абстрактное среднее (логарифмические производные ~
$Lf, ~Lg$, предполагаются неотрицательными).
\Endproc

Важно отметить, что так как функции  из (38) не являются квадратами
средних, то обобщения из теоремы 9  не сводятся к обобщениям из
теоремы 8. Более того,  подынтегральные выражения при выборе (38) в
(33) не являются даже квадратами однородных величин относительно
$f$,\ $g$. Таким образом, мы получаем совершенно неожиданный результат, что интегральный аналог
\textit{необходимой} части известной по дискретному случаю теоремы CDE
Карлица - Дэйкина - Элиезера \hspace{1mm}
\textbf{НЕВЕРЕН!} \ Нами построен класс обобщений интегрального
неравенства Коши--Буняковского такого вида, который запрещён этой
теоремой в дискретном случае для сумм.

Анализ доказательства теоремы CDE выявляет то место, в котором
рассуждения для случая сумм нельзя повторить для случая интегралов.
В теореме CDE используется возможность составить сумму из двух
слагаемых, при этом получается гибридное условие (29), которое мы уже обсуждали. Но для интегралов нет условия,
аналогичному взятию суммы из двух первых слагаемых.

Вместе с тем справедливо следующее утверждение,  дающее тот максимум
информации, который удается получить при попытке хотя бы частично
повторить доказательство необходимой части теоремы CDE. Это удается,
так как для интегралов нетрудно подобрать процедуру, аналогичную
выделению \textit{одного} фиксированного слагаемого из суммы.

\Theorem { 10}
Если справедливо некоторое обобщение
неравенства Коши - Буняковского вида (33), то функции $\Phi_1(f,g),
\Phi_2(f,g)$ связаны в нём тем же соотношением, что и сопряжённые
средние в (11)\,:
$$
\Phi_1(f,g) \cdot \Phi_2(f,g)=f^2(x) \cdot g^2(x).
$$
\Endproc

Однако, как следует из приведенной выше теоремы 9, для интегрального
случая необязательно, чтобы $\Phi_1$ и $\Phi_2$ являлись именно
парой сопряженных \textit{средних}, как в дискретном случае требует
теорема CDE. Например, пара экспоненциальных функций из теоремы 9
вообще не является парой каких--либо средних, так как мы уже
отмечали, что они не являются однородными относительно функций
$f(x), g(x)$, то есть
$$
\Phi_1(\lambda f,\lambda g)\neq \lambda\Phi_1(f,g), \Phi_2(\lambda
f,\lambda g)\neq \lambda\Phi_2(f,g)
$$

В этом
проявляется неожиданное \textit{отличие интегрального случая от
дискретного}.

Выпишем некоторые конкретные оценки. При выборе степенных средних получается

\sc Следствие 9.1. \rm При сделанных предположениях справедливо
следующее обобщение неравенства Коши - Буняковского:
\begin{eqnarray*}
\left(\int^b_af(x)g(x)dx\right)^2\le\int^b_a\exp\left[2\int^x_a
M_\alpha(L
f,L g)\,dt \right]dx\cdot  \\
\cdot\int^b_af^2g^2\exp\left[-2\int^x_a M_\alpha(L f, L
g)\,dt\right]dx\le\int^b_af^2dx\int^b_ag^2dx.
\end{eqnarray*}

При выборе средних Радо получается другое

\sc Следствие 9.2. \rm При сделанных предположениях справедливо
следующее обобщение неравенства Коши - Буняковского:

\begin{eqnarray*}
\left(\int^b_af(x)g(x)dx\right)^2\le\int^b_a\exp\left[2\int^x_a
R_\beta(L
f,L g)\,dt \right]dx\cdot  \\
\cdot\int^b_af^2g^2\exp\left[-2\int^x_a R_\beta(L f, L
g)\,dt\right]dx\le\int^b_af^2dx\int^b_ag^2dx.
\end{eqnarray*}

Аналогично приведенным выше можно выписать вторую серию причудливых
неравенств со средними арифметическим, геометрическим, квадратичным,
логарифмическим, многоэтажно-показательным, а также
арифметико--геометрическим Гаусса и несимметричными весовыми.

При сравнении двух типов построенных нами в теоремах 8 и 9 обобщений
неравенства Коши--Буняковского возникает естественный вопрос об их
взаимозависимости. Представляется, что неравенства  теоремы 9 не
выводятся из неравенств  теоремы 8, так как в них под знаками
интегралов стоят более общие выражения, которые, как отмечалось
выше, не являются однородными относительно функций $f,g$. Можно ли
наоборот вывести неравенства  теоремы 8  из неравенств теоремы 9?

В одном частном случае такой вывод формально возможен. Если в
формуле (38) в качестве среднего выбрать медианту логарифмических
производных, рассматриваемых как дроби, то получится в точности
неравенство другого типа (35) при $\alpha=1$. В случае  выбора обобщённой медианты
получится в точности неравенство (35). Неясно, можно
ли эти рассуждения превратить в строгий вывод одной формулы из
другой, так как медианта --- это не совсем обычное числовое среднее.
Интересно отметить применение при рассмотрении обобщений неравенства
Коши--Буняковского понятия медианты из достаточно далёкой от них
области арифметики, а именно, теории рядов Фарея. Это не случайно,
так как само неравенство Коши--Буняковского находит многочисленные
применения в теории чисел.

\subsec{ 6.3. Основные теоремы об  обобщении интегрального неравенства Коши--Буняковского.}

 В заключении этого раздела приведем наиболее общие формулировки
наших результатов. Это становится возможным с использованием
полученного выше в разделе 4 полного описания средних от двух аргументов.

\Theorem { 11}
Пусть $h(t)$ является произвольной определённой
на всей оси непрерывной чётной функцией, удовлетворяющей следующим
условиям при $t\ge 0$:
$$
h(0)=\frac{1}{2};\ \
\frac{e^{t_1}(e^{t_2}+1)}{e^{t_2}(e^{t_1}+1)}\le
\frac{h(t_1)}{h(t_2)} \le  \frac{(e^{t_2}+1)}{(e^{t_1}+1)},t_1\le
t_2.
$$
Тогда справедливо  обобщение неравенства Коши - Буняковского для
произвольных функций $f(x),\ g(x)$:
$$
\left(\int^b_af(x)g(x)\,dx \right)^2\le \int_a^b{\left[
(x+y)h(\ln\frac{y}{x})\right]}^2\,dx
\cdot \phantom{aaaa} \nonumber\\
$$
$$
\cdot\int_a^bf^2(x)g^2(x)/{\left[(x+y)h(\ln\frac{y}{x})\right]}^2\,dx
\le\int^b_af^2\,dx \int^b_ag^2\,dx.
$$
\Endproc

\Theorem { 12}
Пусть $h(t)$ является произвольной определённой
на всей оси непрерывной чётной функцией, удовлетворяющей следующим
условиям при $t\ge 0$:
\begin{eqnarray*}
h(0)=\frac{1}{2} ,\\
\frac{e^{t_1}(e^{t_2}+1)}{e^{t_2}(e^{t_1}+1)}\ \le \
\frac{h(t_1)}{h(t_2)}& \le & \frac{(e^{t_2}+1)}{(e^{t_1}+1)}, \
t_1\le t_2.
\end{eqnarray*}
Тогда справедливо  обобщение неравенства Коши - Буняковского для
произвольных функций $f(x),g(x)$:
$$
\left(\int^b_af(x)g(x)\,dx \right)^2\le \int_a^b\exp\left(2\int^x_a
(Lf(y)+ Lg(y))h(\ln{\frac{Lf(y)}{Lg(y)}})\,dy\right)\,dx \cdot
$$
$$
\cdot\int_a^b f^2(x)g^2(x)\exp\left(-2\int^x_a \left(Lf(y)+
Lg(y)\right)h(\ln{\frac{Lf(y)}{Lg(y)}})\,dy\right)\,dx \le
$$
$$
\le \int^b_af^2\,dx \int^b_ag^2\,dx.
\phantom{aaaaaaaaaaaaaaaaaaaaaaaaaaa}
$$
\Endproc

 Две последние теоремы носят для данной тематики в определённом смысле окончательный характер. Все рассмотренные выше обобщения интегрального неравенства  Коши--Буняковского получаются из этих теорем как частные случаи
при выборе конкретных функций $h$. Таким образом, теоремы 11 и 12
являются наиболее общими результатами для конструирования средних
частей обобщений рассматриваемого нами вида.

Во всех приведенных выше неравенствах получены условия, при которых
в них достигаются равенства. В некоторых случаях этот вопрос
решается не так просто,
как может показаться.

\subsec{ 6.4. Сравнение обобщений интегрального неравенства Коши--Буняковского.}

Перейдем к  сравнению различных обобщений интегрального неравенства
Коши - Буняковского между собой.

\sc Определение. \rm Введем отношение частичной упорядоченности
на множестве обобщений вида (33): будем записывать $M\prec N$, если
справедливы неравенства
$$
\left(\int^b_afgdx\right)^2\le\int^b_a\left[M(f,g)\right]^2dx\cdot
\int^b_a\left[M^*(f,g)\right]^2dx\le
$$
$$
\le\int^b_a\left[N(f,g)\right]^2dx\cdot
\int^b_a\left[N^*(f,g)\right]^2dx\le\int^b_af^2dx\int^b_ag^2dx,
$$
где $M, N$ --- абстрактные средние, а \ $M^*, N^*$ --- сопряжённые (11) к
ним.

Следует отметить‚ что выполнение введённого отношения
упорядоченности $M\prec N$ в наших терминах означает‚ что обобщение‚
построенное по среднему M‚ лучше приближает левую‚ а обобщение‚
построенное по другому среднему $N$‚
---  правую части интегрального неравенства Коши--Буняковского.

Ясно, что  одним из первых возникает вопрос: упорядочены ли между
собой наши обобщения для средних одного и того же типа, скажем для  $M_\alpha$ или $R_\beta$? Ответы дает

\Theorem { 13}
Пусть $0\le\alpha<\beta\le+\infty$. Тогда
выполнены соотношения
$$M_\alpha\prec M_\beta, R_\alpha\prec R_\beta.$$
\Endproc

 Эта теорема полностью решает поставленный вопрос для  средних Радо или
степенных. Распишем содержащиеся в ней оценки  при $\alpha < \beta$
более подробно:

$$
\left(\int^b_a
fgdx\right)^2\le\int^b_a\left[M_\alpha(f,g)\right]^2dx\cdot
\int^b_a\left[M_{-\alpha}(f,g)\right]^2dx\le
$$
$$
\le
\int^b_a\left[M_\beta(f,g)\right]^2dx\cdot\int^b_a
\left[M_{-\beta}(f,g)\right]^2dx\le\int^b_af^2dx\int^b_a
g^2dx,\eqno(39)
$$
$$
\left(\int^b_afgdx\right)^2\le\int^b_a\left[R_\alpha(f,g)
\right]^2dx\cdot\int^b_a\left[R_{\alpha}^*(f,g)\right]^2dx\le
$$
$$
\le\int^b_a\left[R_\beta(f,g)\right]^2dx\cdot\int^b_a
\left[R_{\beta}^*(f,g)\right]^2dx\le\int^b_af^2dx\int^b_ag^2dx.
$$

В частности‚ из (39) следует‚ что нами получено бесконечное число
неравенств‚ которые обобщают неравенство Милна с $(M_2, M_{-2})$ как
с точки зрения лучшего приближения к скалярному произведению (при
$0<\alpha<2$)‚ так и к правой части (при $2<\alpha<+\infty$).

Результаты по сравнению обобщений со средними из различных шкал вида
$M_{\alpha}\prec R_{\beta}$‚ или тем более вида
$M_{\alpha}\prec\mu,\  R_{\beta}\prec\mu$‚ (а также обратные к ним)
получены автором  только для некоторых частных случаев, общий способ
доказательства таких неравенств мне неизвестен.

Перейдем к  более содержательному вопросу о сравнении обобщений из
теоремы 9 с экспоненциальными выражениями. Для них также будем
употреблять символ отношения частичной упорядоченности $M\prec N$  в
прежнем смысле.

\Theorem { 14}
Справедливы следующие отношения:\\
$$
M_1 \prec M_{\alpha} \prec M_{\beta}, 1\le\alpha\le\beta\le+\infty,
$$
$$
M_1\prec M_{\beta}\prec M_{\alpha},
-\infty\le\alpha\le\beta\le 1.
$$
\Endproc

\medskip  Аналогично решен вопрос и для средних Радо.

Из теоремы 14 следует, что остаётся изучить сравнения обобщений со
степенными средними, один из порядков которых больше, а другой
меньше единицы. Для этого случая вначале были  получены некоторые
отдельные соотношения:
$$
M_2\prec M_0, \   M_3\prec M_{-1}, \ M_{\frac{3}{2}}\prec
M_{\frac{1}{2}},
$$
а также ряд других. В частности, усиления со средними
$M_{\frac{1}{2}}$ и $M_2$ оказались не сравнимы. Такая возможность
является достаточно неожиданной и, по мнению автора, ещё с одной
стороны подчёркивает содержательность рассматриваемого здесь круга
задач.

Затем для общего случая автором совместно с А.~Кореновским была
доказана

\Theorem { 15}
Пусть  $\alpha<1, \beta>1, \alpha +
\beta
>2$. Тогда усиления со средними $M_\alpha$ и $M_\beta$ не
сравнимы.
Пусть  $\alpha<1, \beta>1, \alpha + \beta \le 2$. Тогда $M_\beta
\prec M_\alpha$.
\Endproc

В результате все доказанные отношения или их опровержения можно
собрать на специальной диаграмме достаточно внушительных размеров. В
частности, выделим такие цепочки:
\begin{eqnarray*}
 R_1=M_1\prec
M_{\frac{3\ln 2}{2\ln\frac{5}{2}}}\prec R_{\frac{3}{2}}\prec
M_{\frac{7}{6}}\prec M_{\frac{2\ln 2}{\ln 3}}\prec R_2 \prec
M_{\frac{4}{3}} \ldots , \\
R_1=M_1\prec R_{\frac{1}{2}}\prec M_{\frac{5}{6}}\prec
M_{\ln_2}\prec R_0 \prec M_{\frac{2}{3}}\prec M_{\frac{1}{2}} = \smallskip \\
= R_{\frac{1}{2}}\prec R_{-1}\prec \mu(M_1, M_0)\prec
R_{-2}=M_0\prec R_{-5} ...
\end{eqnarray*}

Отношения с итерационным средним $\mu$ получены с использованием
неравенства из пункта 4.3.

Фактически надо признать, что  пока не удалось найти достаточно
общий метод доказательства подобных теорем, и каждое неравенство
 доказывается своим способом. В идеале было бы найти такой метод,
чтобы из каждого
неравенства между двумя средними сразу выводились отношения порядка
между соответствующими обобщениями неравенства Коши--Буняковского, или делался вывод об их несравнимости.

\section{7. Приложения обобщений дискретных неравенств Коши--Буняковского.}

\subsec{ 7.1 Неравенства Коши--Буняковского для $q$--интеграла Джексона.}
Следующий пример относится к области коммутативного $q$--Анализа. Это такой раздел Анализа, в котором есть свои понятия $q$--производной, $q$--интеграла и соответствующие обобщения специальных функций, которые при $q\to 1$ переходят в обычные определения~[7]. Так, например, $q$--интеграл Джексона вводится по формуле
$$
\int_0^1 f(t)\,d_q t=(1-q)\sum_{k=0}^\infty f(q^k) q^k.
$$
Этот интеграл на самом деле не интеграл, а ряд. Для него получается
\Theorem { 16}
 Пусть $M$ - произвольное абстрактное среднее‚ $M^*$ -
сопряженное к нему. Тогда справедливо обобщение неравенства Коши--Буняковского в терминах $q$ - интеграла Джексона:
$$
\left(\int_0^1f(t)g(t)d_qt\right)^2\le \left(\int_0^1(M(f(t),
g(t)))^2d_qt\right)\cdot
$$
$$
\cdot\left(\int_0^1(M^*(f(t),
g(t)))^2d_qt\right)\le
\left(\int_0^1f^2(t)d_qt\right)^2\left(\int_0^1g^2(t)d_qt\right)^2.
\phantom{xxxxxxxxxxxxxxxxx}
$$
\Endproc

\subsec{ 7.2. Неравенства Коши--Буняковского в пространствах Лоренца.}

Сравнительно малоизвестным является тот факт, что в пространствах
Лоренца~[3] также выполняется неравенство Коши--Буняковского, но в
обратную сторону. Замечательное доказательство Бохнера основано
на интегральном представлении полунормы по световому конусу. Данное
неравенство является частным случаем известных неравенств Ацеля. Именно на данном неравенстве основано математическое
обоснование так называемого "парадокса близнецов"\  в специальной
теории относительности~[3].

Для этого случая получена

\Theorem { 17}
Пусть $A(x, y)$ есть средняя часть
произвольного обобщения дискретного неравенства Коши--Буняковского
вида
$$
\left(\sum_{k=1}^nx_k\cdot y_k\right)^2\le A(x,
y)\le\left(\sum_{k=1}^nx_k^2\right)\cdot\left(\sum_{k=1}^ny_k^2\right),
$$
и выполнены условия
$$
x_0^2-\sum_{k=1}^nx_k^2\geq 0,\  y_0^2-\sum_{k=1}^ny_k^2\geq 0.
$$
Тогда справедливо следующее уточнение дискретного неравенства
Коши--Буняковского в пространствах Лоренца:
$$
\left(x_0y_0-\sum_{k=1}^nx_ky_k\right)^2\geq \left(x_0y_0-\sqrt{A(x,
y)}\right)^2 \ge
$$
$$
\ge\left(x_0^2-\sum_{k=1}^nx_k^2\right)\left(y_0^2-\sum_{k=1}^ny_k^2\right).
$$
\Endproc

 Аналогично обычному случаю из последнего неравенства выводится
обобщение неравенства Минковского для пространств Лоренца. Его же
можно применить для дальнейшего уточнения оценок в известных
геометрических задачах~[70,98]. Вопрос о справедливости
подобных неравенств  в пространствах Понтрягина с произвольной
сигнатурой (например, $+ + - -$), насколько известно автору, не
изучен.

\subsec{ 7.4. Принцип неопределённости для ДПФ.}

Рассмотрим дискретное преобразование Фурье (ДПФ), которое определяется на $n$--мерных векторах по формулам
$$
b_j=\frac{1}{\sqrt n}\sum_{k=0}^{n-1}a_k w^{-jk},
a_m=\frac{1}{\sqrt n}\sum_{k=0}^{n-1}b_k w^{mk},
$$
где $w=\exp(\frac{2\pi i}{n})$--первообразный корень из единицы.

Сейчас не будем агитировать, что это одно из самых важных понятий в математике, а приведём такую форму принципа неопределённости, которая доказана в~[74]. Доказательство там не очень длинное, но я бы мягко говоря не назвал его понятным. Случай равенства указан один из возможных, но не все. Дадим своё полное доказательство.

\Theorem { 18} Пусть $A,B$--число ненулевых компонент в соответствующих векторах, $n$--их размерность. Тогда справедливо соотношение неопределённости
$AB\ge n$, причём равенство достигается тогда и только тогда, когда в одном из векторов все ненулевые компоненты равны.
\Endproc

\beginproof~Получаем, последовательно применяя неравенства треугольника и Коши--Буняковского
$$
|b_j|\le \frac{1}{\sqrt n}\sum_{a_k\neq0}|a_k|\le \frac{1}{\sqrt n}\sqrt A
\sqrt{\sum_{a_k\neq0}|a_k|^2}\le \sqrt{\frac{A}{n}}\|a\|_2.
$$
Далее возводим в квадрат, складываем ненулевые компоненты, применяем равенство Парсеваля--Стеклова---всё!!! Случай равенства  следует из условий равенства для обычного неравенства Коши--Буняковского.
~\endproof

Смысл этого утверждения в следующем: последовательность и её ДПФ не могут иметь одновременно слишком много нулей.

Теперь можно выписывать различные обобщения приведённого неравенства по разработанной методике.

В заключение этого раздела отметим, что
дискретном случае для абстрактных средних, удовлетворяющих набору
свойств 1) -- 4) из определения, можно доказать аналог наиболее
общей теоремы 11 с произвольной функцией h. А вот аналог теоремы 9
пока не найден.
По моему мнению, это очень трудная задача, требующая нахождения
новой идеи. Но она имеет принципиальное значение для рассматриваемой
здесь методики, так как позволит не только констатировать различие
дискретного и интегрального случаев, но и построить более богатый
набор обобщений дискретного неравенства Коши--Буняковского, выйдя
за рамки ограничений теоремы CDE.

Разумеется, любое уточнение неравенств Коши - Буняковского в
дискретном или интегральном случаях ведет к уточнению
соответствующих неравенств Минковского.

\section{8. Приложения обобщений интегральных неравенств Коши--Буняковского.}

В этом пункте мы ограничимся приложениями обобщений интегральных неравенств Коши--Буняковского к выводу оценок для специальных функций. Случай гамма--функций Эйлера, неполных гамма--функций, некоторых функций Ломмеля и Миттаг--Лефлера рассмотрен в~[37--39]. Здесь рассмотрим вывод неравенств для полных эллиптических интегралов Лежандра первого рода, неравенства для которых автор начал изучать в~[44]. Отметим, что в указанной работе одной из первых была установлена связь между гипергеометрическими функциями и полилогарифмами; сейчас эта тематика активно разрабатывается.

В~[37--39] приводятся общие теоремы,  позволяющие каждое из
полученных ранее обобщений неравенства Коши-Буняковского развить в
некоторую итерационную процедуру,  генерирующую бесконечное число
новых последовательных обобщений того же неравенства, каждое
следующее из которых в определённом смысле точнее предыдущего.

При этом  основную роль начинают играть именно левые
очевидные (но не тривиальные!) части наших обобщений (33).
Основная идея проста: раз неравенство Коши--Буняковского позволяет
уточнять само себя, то этот приём можно просто повторять сколько
угодно раз. Как показано ниже, на этом пути удаётся, например,
получать эффективные оценки специальных функций, заданных своими
интегральными представлениями.

В качестве приложения одной из простейших подобных итерационных процедур уточнения неравенства Коши--Буняковского, полученных в~[37--39],  рассмотрим  неравенства для полных
эллиптических интегралов Лежандра первого рода~[2,63]. Эта
гипергеометрическая функция имеет интегральное представление при
$0\le x< 1$
\begin{equation*}
K(x)=\int_0^1\frac{dt}{\sqrt{(1-t^2)(1-x^2t^2)}}=\frac{\pi}{2}\cdot
{}_2F_1\left(\left.
\begin{array}{cc}
  1/2 & 1/2 \\
  1 & {} \\
\end{array}
\right| x^2 \right)\eqno(40)
\end{equation*}
и особенность с логарифмической  асимптотикой при $x\rightarrow 1$.

Неравенства для эллиптических интегралов различных типов --- это
целый раздел теории специальных функций, насчитывающий несколько
сотен работ. Достаточно подробный перечень приведён в~[37--39], здесь же мы процитируем только работы, написанные при участии автора~[44,80--86].

Для использования указанной выше итерационной
процедуры  ключевую роль играет правильное "расчленение" \
интеграла  в скалярное произведение. Осуществим следующий
выбор:
$$
K^2(x)=\left(\int_0^1f(t)g(t)dt\right)^2,
$$
$$
f=\frac{1}{(1+t)^{1/2}(xt^2-(x+1)t+1)^{1/4}},
$$
$$
g=\frac{1}{(1+xt)^{1/2}(xt^2-(x+1)t+1)^{1/4}}.
$$
Он обусловлен необходимостью: во--первых, предусмотреть сходимость
интегралов на всех шагах итерационной процедуры; во--вторых,
выполнить равенство $f=g$ при $x=1$, которое позволяет надеяться на
повышенную точность оценок при $x\rightarrow 1$; в--третьих,  явно
вычислять получаемые интегралы.

 Мы приведем
список вычисленных оценок для упомянутой итерационной процедуры  при использовании в качестве уточняющего неравенство Милна. При этом нетривиальные оценки получаются уже из самого
неравенства Коши - Буняковского на начальном шаге  за счет удачного выбора
 "расчленяющих" \  функций.
\begin{eqnarray*}
G_0(x)=\frac{1}{\sqrt{(x+1)2x}}\ln\left(\frac{2\sqrt{2x(x+1)}+3x+1}{1-x}\right),
\\
L_0(x)=\frac{1}{\sqrt{2(x+1)}}\ln\left(\frac{2\sqrt{2(x+1)}+x+3}{1-x}\right),
\\
G_1(x)=\frac{1}{2}\left[\frac{1}{\sqrt{2x(x+1)}}
\ln\left(\frac{2\sqrt{2x(x+1)}+3x+1}{1-x}\right)\right]+
\\
+\frac{1}{\sqrt{2(x+1)}}\ln\left(\frac{2\sqrt{2(x+1)}+x+3}{1-x}\right),
\\
L_1(x)=\frac{2}{\sqrt{(x+3)(3x+1)}}\ln\left(\frac{\sqrt{(x+3)(3x+1)}+2x+2}{1-x}\right),
\\
G_2(x)=\frac{\left(\frac{\sqrt{5}+1}{2\sqrt{5}}\right)} { \sqrt{
(x+4+\sqrt{5}) \left( \left(\frac{7+\sqrt{5}}{2}
\right)x+\frac{3+\sqrt{5}}{2}\right) }
} \times\\
\times\ln\left(\frac{2\sqrt{(x+4+\sqrt{5})\left(\left(\frac{7+\sqrt{5}}
{2}\right)x+\frac{3+\sqrt{5}}{2}\right)}
+\left(\frac{9+\sqrt{5}}{2}\right)x+\frac{11+3\sqrt{5}}{2}}{\left(\frac{5+\sqrt{5}}
{2}\right)(1-x)}\right)+\\
+\frac{ \left( \frac{\sqrt{5}-1}{2\sqrt{5}} \right) } { \sqrt{
(x+4-\sqrt{5}) \left( \left( \frac{7-\sqrt{5}} {2} \right)
x+\frac{3-\sqrt{5}}{2} \right) } }\times \\
\times\ln\left(\frac{2\sqrt{(x+4-\sqrt{5})\left(\left(\frac{7-\sqrt{5}}
{2}\right)x+\frac{3-\sqrt{5}}{2}\right)}
+\left(\frac{9-\sqrt{5}}{2}\right)x+\frac{11-3\sqrt{5}}{2}}{\left(\frac{5-\sqrt{5}}
{2}\right)(1-x)}\right),
\end{eqnarray*}
\begin{eqnarray*}
L_2(x)=\frac{1}{\sqrt{2x+2}}\ln\left(\frac{2\sqrt{2x+2}+x+3}{1-x}\right)+ \\
+\frac{1}{\sqrt{2x(x+1)}}\ln\left(\frac{2\sqrt{2x(x+1)}+3x+1}{1-x}\right)+ \\
+\frac{1}{\sqrt{(x+3)(3x+1)}}\ln\left(\frac{\sqrt{(x+3)(3x+1)}+2x+2}{1-x}\right).
\end{eqnarray*}

\Theorem { 18}
Справедливы  следующие двусторонние неравенства для
полного эллиптического интеграла Лежандра:
$$
L_0(x) \le L_1(x) \le L_2(x) \le K(x) \le G_2(x) \le G_1(x) \le G_0(x).
$$
\Endproc

Расчеты показывают высокую точность полученных двусторонних
неравенств. Они также дают всё больше истинных слагаемых из асимптотики при $x\to 1$. Вычисления можно продолжать и дальше. Было бы интересно разобраться в закономерностях образования  всё более точных, но громоздких оценок.

Приведенные  результаты для эллиптических интегралов Лежандра могут
быть применены к различным оценкам с этими функциями. В результате
получим оценки уже через элементарные функции, но при этом без потери
асимптотики. Таким образом могут быть модифицированы многие
неравенства с $K(x)$: теоремы площадей, формулы
остаточного члена интерполяционной формулы Эрмита, результаты
для однолистных функций и проблемы коэффициентов,
оценки сингулярных интегралов на единичной окружности и многочленов
Фабера, остатков рядов Тэйлора и их констант Лебега в
комплексной области, некоторые оценки в теории рациональной
аппроксимации, формула Стечкина для Фурье - уклонения класса
функций Соболева от класса тригонометрических полиномов.
Важной областью исследований является использование оценок для
эллиптических интегралов и более общих гипергеометрических функций
нескольких переменных, таких как функции Аппеля и Лауричеллы, в
теории потенциала и особенно в задачах, связанных с ёмкостями. Подробные ссылки см. в [37--39].

При некоторых формулировках плоской теории упругости было замечено,
что ограничения на основные тензоры принимают в точности вид
неравенств Коши - Буняковского~[39]. Таким образом, с
использованием наших обобщений эти ограничения могут быть несколько
ослаблены.

Для интегрального оператора $\int_a^b K(x,y)f(y)\,dy$ неравенство Коши--Буняковского даёт естественный способ оценки норм в пространстве $L_2$. Полученные результаты позволяют выписать многочисленные обобщения этих оценок. Интересным представляется применить данный способ к конкретным интегральным операторам, а именно, операторам преобразования (ОП). Подробная теория ОП изложена в предыдущем обзоре автора для этого сборника 2007 г.~[45], см. также~[10--13,22,46--50]. Другая возможность заключается в применении полученных результатов к интегральным операторам, рассматриваемым в более общих ситуациях, например, как в~[14,20,31,59,64].

Многие результаты в теории вероятностей и математической статистике имеют форму неравенств, которые по существу являются именно неравенствами Коши--Буняковского в специальных терминах. Например, такую форму имеют неравенства для коэффициентов корреляции или регрессии. Эти результаты могут быть также уточнены по схеме настоящей работы для дискретных или непрерывных случайных величин.

В заключение этого пункта отметим, что основные результаты данной работы
докладывались автором на различных
международных и российских конференциях, начиная с 1996 года~[106,8,51--57].

Автор выражает благодарность за обсуждения результатов работы А.\,Кореновскому (Одесса), Л.\,Минину
(Воронеж), А.\,Лободе (Воронеж), И.\,Половинкину (Воронеж), В.\,Родину (Воронеж), А.\,Боровских (Москва), Д.\,Карпу (Владивосток), В.\,Фетисову (Шахты), С.\,Калинину (Ижевск), Т.\,Поганжу и А.\,Баричу (Хорватия) Л.\,Брылевской (Санкт-Петербург) за
помощь в копировании оригинальной работы В.Я.\,Буняковского (к
сожалению, рукописи не только горят, но и рассыпаются от ветхости,
что и произошло в данном случае при копировании). Мне также не забыть  обсуждения
возможных приложений данной тематики с проф. В.\,В. Сысоевым.
Автор благодарит Анатолия Георгиевича Кусраева за внимание к моей работе и предложение написать этот обзор, а также обсуждение возможности перенесения некоторых результатов на случай векторных решёток.

\par\bigskip\centerline{\bf Литература}\smallskip

 \begin{enumerate}

 \itemsep=0pt\parskip=0pt

\bib{Анциферова~Г.~А., Ситник~С.~М.}{Некоторые обобщения неравенства Юнга~/\!/Вестник Воронежского института МВД России. Воронеж, 1999.---\No~2(4).---С.~161--164.}
\bib{Бейтмен~Г., Эрдейи~А.}{Высшие трансцендентные функции. Т.1--3.---М., Наука.}
\bib{Беккенбах~Э., Беллман.~Р.}{Неравенства.---М.: Мир,
1965 (1 изд 1961 г.).---276~c.}
\bib{Блох~А.~Ш., Неверов~Г.~С.}{Решение неравенств.---Минск: Изд-во М-ва высш., сред. спец. и проф. образования БССР, 1962.---43~с.}
\bib{Борвейн~Дж.~М., Борвейн~П.~Б.}{Рамануджан и число  пи~/\!/В мире науки.  Издание на русском языке.---1988.---\No~4.---C.~58--66.}
\bib{Буренков~В.~И.}{Функциональные пространства. Основные интегральные неравенства, связанные с пространствами \(L_p\).---М.: УМН, 1989.---96~с.}
\bib{Гаспер~Дж., Рахман~М.}{Базисные гипергеометрические ряды.---М.: Мир, 1993.---327~c.}
\bib{Дикарева~Е.~В., Ситник~С.~М.}{Об одном обобщении неравенства Коши---Буняковского с максимумом и минимумом.~/\!/~В сб.:~Труды участников международной школы -семинара
по геометрии и анализу памяти Н.В.Ефимова. Абрау-Дюрсо.---Ростов-на-Дону, Южный Федеральный университет, 2008.---С.~107--108.}

\bib{Калинин~С.~И.}{Средние величины степенного типа. Неравенства Коши и Ки Фана.---Киров, 2002.---368~с.}
\bib{Катрахов~В.~В., Ситник~С.~М.}{Краевая задача для стационарного уравнения Шрёдингера с сингулярным потенциалом~/\!/~ДАН СССР.---1984.---т.278.---№4.---С.797--799.}
\bib{Катрахов~В.~В., Ситник~С.~М.}{Оценки решений Йоста одномерного уравнения Шредингера с сингулярным потенциалом~/\!/~ДАН СССР.---1995.---Т.~340.---№~1.---С.~18--20.}
\bib{Катрахов~В.~В., Ситник~С.~М.}{Метод факторизации в теории операторов преобразования~/\!/~В сб.:~(Мемориальный сборник памяти Бориса Алексеевича Бубнова). Неклассические уравнения и уравнения смешанного типа.
(ответственный редактор В.~Н.~Врагов).---1990,~Новосибирск.---С.~104--122.}
\bib{Катрахов~В.~В., Ситник~С.~М.}{Композиционный метод
построения В--эллиптических, В--гиперболических и В--параболических операторов преобразования~/\!/~ДАН СССР.---1994.---Т.~337.---№~3.---С.~307--311.}
\bib{Коробейник~Ю.~Ф.}{О разрешимости в комплексной области некоторых общих классов линейных интегральных уравнений.---Ростов--на--Дону:~2005.---245~с.}
\bib{Коровкин~П.~П.}{Неравенства.---М.: Наука, 1966 (3 изд.).---56~с.}
\bib{Крыжановский~Д.~А.}{Элементы теории неравенств.---М.-Л.: Объединенное научно-техническое изд. НКТП СССР (ОНТИ), 1936.---112~c.}
\bib{Кусраев~А.~Г.}{Векторная двойственность и ее приложения.---Новосибирск: Наука, 1985.---256~с.}
\bib{Кусраев~А.~Г., Кутателадзе~С.~С.}{Субдифференциальное исчисление.---Новосибирск, Наука, СО, 1987.---224~c.}
\bib{Кусраев~А.~Г., Кутателадзе~С.~С.}{Субдифференциалы. Теория и приложения. Части 1 и 2.---Новосибирск, ИМ им. С.~Л.~Соболева, 2002, 2003.}
\bib{Кусраев~А.~Г.}{Мажорируемые операторы.---М.:~Наука,~2003.---619~с.}
\bib{Кутателадзе~С.~С., Рубинов~А.~М.}{Двойственность Минковского и её приложения.---М.: Наука, 1976.---257~c.}

\bib{Ляховецкий~Г.~В., Ситник~С.~М.}{Операторы преобразования
Векуа--Эрдейи--Лаундеса~/\!/~Препринт института автоматики и
процессов управления Дальневосточного отделения РАН.---Владивосток:~ДВО РАН,
1994.---24~с.}
\bib{Маркус~М., Минк~Х.}{Обзор по теории матриц и матричных
неравенств.---М.: Наука, 1972.---232~с.}
\bib{Маршалл~А., Олкин~И.}{Неравенства: теория мажоризации и ее
приложения.---М.: Мир, 1983.---575~c.}
\bib {Маслов~В.~П.}{Нелинейное среднее в экономике~/\!/Математические заметки.---2005.---Вып.~3., Т.~78.---С.~377--395.}
\bib {Маслов~В.~П.}{Нелинейное финансовое осреднение, эволюционный процесс и законы эконофизики~/\!/Теория вероятностей и приложения.---2004.---\No~49:2.---С.~269–-296.}
\bib {Маслов~В.~П.}{О нелинейности осреднений в финансовой математике~/\!/Математические заметки.---2003.---Вып.~6., Т.~74.---С.~944–-947.}
\bib {Маслов~В.~П.}{Аксиомы нелинейного осреднения в финансовой математике и динамика курса акций~/\!/Теория вероятностей и приложения.---2003.---\No~48:4.---С.~800–-810.}
\bib{Натансон~И.~П.}{Теория функции вещественной переменной.---М.: Наука, 1974.---480~с.}
\bib{Невяжский~Г.~Л.}{Неравенства: Пособие для учителей.---М.: Учпедгиз,
 1947,---204~c.}

\bib{Пасенчук~А.~Э.}{Абстрактные сингулярные операторы.---Новочеркасск,~1993. }
\bib{Полиа~Г., Сеге~Г.}{Изопериметрические неравенства в математической физике.--- М.: ГИФМЛ, 1962.---336~с.}
\bib{Половинкин~Е.~С., Балашов~М.~В.}{Элементы выпуклого и сильно выпуклого анализа.---М.: ФМЛ, 2004.---416~с.}
\bib{Рокафеллар~Р.}{Выпуклый анализ.---М.: Мир, 1973.---470~с.}
\bib{Седракян~Н.~М., Авоян~А.~М.}{Неравенства. Методы доказательства.---М.: Физматлит, 2002.---256~с.}
\bib{Сивашинский~И.~Х.}{Неравенства в задачах.---М.: Наука, 1967.---304~с.}
\bib {Ситник~С.~М.}{Уточнение интегрального неравенства Коши--Буняковского~/\!/Вестник Самарского гос. тех. университета. Сер.
"Физико-математические науки"\.---2000.---\No~9.---С.~37--45.}
\bib {Ситник~С.~М.}{Некоторые приложения уточнений неравенства
Коши-Буняковского~/\!/Вестник Самарской
государственной экономической академии.---2002.---\No~1(8).---С.~302--313.}
\bib{Ситник~С.~М.}{Обобщения неравенств Коши-Буняковского методом средних значений и их приложения~/\!/~Чернозёмный альманах научных исследований. Серия "Фундаментальная математика".---2005---№1~(1).---С.~3--42.}
\bib {Ситник~С.~М.}{Сколько неравенств заключено в неравенстве Юнга?~/\!/Труды Всероссийской заочной научно--практической конференции
"Актуальные проблемы обучения математике"\,,
посвящённой 155--летию со дня рождения
Андрея Петровича Киселёва.---Орёл, Орловский государственный университет,
2007.---С.~464--469.}
\bib{Ситник~С.~М.}{Теоремы Радо и среднее логарифмическое.~/\!/Препринт института автоматики и процессов управления
Дальневосточного отделения РАН.---Владивосток, 1992.---16~c.}
\bib{Ситник~С.~М.}{Неравенства для среднего логарифмического и итерационных средних~/\!/Препринт института автоматики и процессов управления Дальневосточного отделения РАН.---Владивосток, 1992.---14~c.}
\bib{Ситник~С.~М.}{Обобщения неравенства Коши--Буняковского в пространствах
с индефинитной метрикой~/\!/~В сб.:~Материалы шестой Казанской международной летней школы-конференции. Теория функций, её приложения и смежные вопросы.---Труды математического центра имени Н.~И.~Лобачевского.---Том 19.---Казань,~2003.---C.~202--203.}
\bib{Ситник~С.~М.}{Неравенства для полных эллиптических интегралов Лежандра~/\!/~~Препринт ИАПУ ДВО РАН.---Владивосток,~1994.---17~С.}
\bib{Ситник~С.~М.}{Операторы преобразования и их приложения.---Исследования по современному анализу и математическому моделированию.
Отв. ред. Коробейник~Ю.~Ф., Кусраев~А.~Г.
Владикавказ: Владикавказский научный центр РАН и РСО--А.
2008.---C.~226--293.}
\bib{Ситник~С.~М.}{Метод факторизации операторов преобразования
в теории дифференциальных уравнений~/\!/~Вестник Самарского Государственного Университета (СамГУ) — Естественнонаучная серия.---2008.---\No~8/1 (67).---С.~237--248.}
\bib{Ситник~С.~М.}{Унитарность и ограниченность операторов Бушмана--Эрдейи нулевого порядка гладкости~/\!/~Препринт ИАПУ ДВО РАН.---Владивосток,~1990.---45~С.}
\bib{Ситник~С.~М.}{Факторизация и оценки норм  в весовых лебеговых пространствах операторов Бушмана-Эрдейи~/\!/~ДАН СССР.---1991.---т.320.---№6.---С.1326--1330.}
\bib{Ситник~С.~М.}{Оператор преобразования  и  представление Йоста для уравнения с сингулярным потенциалом~/\!/~~Препринт ИАПУ ДВО РАН.---Владивосток,~1993.---21~С.}
\bib{Ситник~С.~М.}{Построение операторов преобразования Векуа--Эрдейи--Лаундеса~/\!/~В сб.: Тезисы докладов международной конференции 'Дифференциальные уравнения, теория функций и приложения', посвящённой 100--летию со дня рождения академика И.~Н.~Векуа.---Новосибирск:~НГУ,~2007.---С.~469--470.}
\bib{Ситник~С.~М.}{Метод получения  последовательных уточнений неравенства Коши--Буняковского и его применения к оценкам специальных функций~/\!/~В сб.:~Современные методы теории краевых задач.
Материалы Воронежской  весенней математической школы 'Понтрягинские чтения-VII'.---Воронеж,~1996.---С.~164.}
\bib{Ситник~С.~М.}{Обобщения неравенств Коши-Буняковского и их приложения~/\!/~В сб.:~Abstracts of Bunyakovsky International Conference.---Kyiv,~2004.---C.~119--120.}
\bib{Ситник~С.~М.}{Об уточнениях интегрального неравенства Коши-Буняковского~/\!/~В сб.:~Современные проблемы теории функций и их приложения. Тезисы докладов 11--й Саратовской зимней школы, посвящённой памяти выдающихся профессоров МГУ Н.~К.~Бари  и  Д.~Е.~Меньшова.---Саратов,~2002.---C.~ 191--192.}
\bib{Ситник~С.~М.}{Обобщения неравенств Коши-Буняковского и их приложения~/\!/~В сб.:~Abstracts of Bunyakovsky International Conference.---Kyiv,~2004.---C.~119--120.}
\bib{Ситник~С.~М.}{Уточнения интегрального неравенства Коши-Буняковского и их
приложения к дифференциальным уравнениям~/\!/~В сб.:~ Международная конференция 'Дифференциальные уравнения и смежные вопросы', посвящённая памяти И.~Г.~Петровского. Тезисы докладов.---Москва,~МГУ.---2004.---С.~210.}
\bib{Ситник~С.~М.}{О некоторых обобщениях неравенства Коши - Буняковского~/\!/~В сб.:~" Международная конференция 'Анализ и особенности', посвящённая 70--летию В.И. Арнольда ".
Тезисы докладов.---Москва, МИАН.---2007.---С.~92--94.}
\bib{Ситник~С.~М.}{Обобщения неравенства Коши --- Буняковского с
использованием средних~/\!/~В сб.:~ Международная конференция 'Дифференциальные уравнения и смежные вопросы', посвящённая памяти И.~Г.~Петровского. Тезисы докладов.---Москва,~МГУ.---2007.---С.~297--298.}

\bib{Тихомиров~В.~М.}{Выпуклый анализ.---М.: ВИНИТИ, 1987.---С.~5--101.}
\bib{Фетисов~В.~Г.}{Операторы и уравнения в локально ограниченных пространствах~/\!/В книге:
Операторы и уравнения в линейных топологических пространствах.---Владикавказ:
 Изд-во ВНЦ РАН,~2006.---С.~7--142.}
\bib{Харди~Г.~Г., Литтлвуд~Дж.~Е, Полиа~Г.}{Неравенства.---М.: ИЛ,
 1948.---456~с.}

\bib {Acz\'{e}l~J.,  Dar\'{o}czy~Z.}{Measures of information and their characterizations.---Academic Press, 1975.---234~p.}
\bib{Anderson~G.D, Vamanamurthy~M.~K., Vuorinen~M.~K.}{Conformal Invariants, Inequalities, and Quasiconformal Maps.---Wiley, 1997.---536~p.}
\bib{Andrews~G.~E., Askey~R., Roy~R.}{Special functions.---Cambridge University Press, 1999,---681~p.}
\bib{Appell J.M., Kalitvin A.S., Zabrejko P.P.
}{Partial Integral Operators and Integro-Differential Equations.---N.~Y.: Marcel Dekker, 2000.---560~p.}
\bib{Berggren~B., Borwein~J., Borwein~P.}{Pi: A Source Book.---Springer, 1967.---718~p.}
\bib{Borwein~J.~M., Borwein~P.~B.}{Pi and the AGM.---Wiley,
1987.---432~p.}
\bib{Borwein~P.}{Computational Excursions in Analysis and Number Theory.---Springer, 2002.---272~p.}
\bib{Borwein~J.~M.,  Bailey~D.~H.}{Mathematics by Experiment: Plausible reasoning in the 21st century.---A.K. Peters Ltd, 2008.---393~p.}
\bib{Borwein~J.~M.,  Bailey~D.~H.}{Experimentation in Mathematics: Computational paths to discovery.---A.K. Peters Ltd, 2004.---368~p.}
\bib{Bottema~O., Djordjevic~R.~E., Mitrinovi\'{c}~D.~S., Vasi\'{c}~P.~M.}{Geometric Inequalities.---Groningen, 1969.---151~p.}
\bib{Bullen~P.~S., Mitrinovi\'{c}~D.~S., Vasi\'{c}~P.~M.}{Means and Their Inequalities.---D.Reidel Publishing Company, Dordrecht, 1988.---480~p.}
\bib{Bullen~P.~S.}{Handbook of Means and Their Inequalities.---Kluwer, 2003.---587~p.}
\bib{Bullen~P.~S.}{A dictionary of inequalities. (Pitman
Monographs and Surveys in Pure and Applied Mathematics).---CRC Press, 1998, Vol. 97.---283~p.}
\bib{Donoho~D.~L., Huo~X.}{Uncertainty Principles and Ideal Atomic Decomposition~/\!/IEEE Transactions on Information Theory.---2001.---Vol.~47, \No~7.---P.~2845--2862.}

\bib{Dragomir~S.~S.}{A Survey on Cauchy--Buniakowsky--Schwartz Type
Discrete Inequalities.---RGMIA monographs, 2003.---214~p.}
\bib{Dragomir~S.~S.}{Advances in Inequalities of the Schwarz, Gruss and
Bessel Type in Inner Product Spaces.---RGMIA monographs, 2004.---257~p.}
\bib{Dragomir~S.~S.}{Advances in Inequalities of the Schwarz, Triangle and
Heisenberg Type in Inner Product Spaces.---RGMIA monographs, 2004.---283~p.}
\bib{Dragomir~S.~S., Pearce~C.~E.~M.}{Selected Topics on Hermite--Hadamard Inequalities and Applications.---RGMIA monographs, Victoria University, 2002.---361~p.}

\bib{Fumio~H., Kosaki~H.}{Means of Hilbert Space Operators.---Springer Lecture Notes in Mathematics, 2003.---\No~1820.---148~p.}
\bib{Karp~D., Savenkova~A., Sitnik~S.~M.}{Series expansions for the third incomplete elliptic integral via partial fraction decompositions~/\!/Journal of Computational and Applied Mathematics.---2007.---Vol.~207.---No.~2.---P.~331--337.}
\bib{Karp~D., Sitnik~S.~M.}{Asymptotic approximations for the first incomplete elliptic integral near logarithmic singularity~/\!/~Journal of Computational and Applied Mathematics.---2007.---Vol.~205.---P.~186--206.}
\bib{Karp~D., Sitnik~S.~M.}{Inequalities and monotonicity of ratios for generalized hypergeometric function~/\!/~J. Approx. Theory.---2009.---Vol.~161.---14~p.}
\bib{Karp~D., Sitnik~S.~M.}{Log-convexity and log-concavity of hypergeometric-like functions~/\!/J. Math. Anal. Appl.---2009 (принята к печати).---14~p.}
\bib{Karp~D., Sitnik~S.~M.}{Two-sided inequalities for generalized hypergeometric function~/\!/~RGMIA Research Report Collection.---2007.--- 10(2).---Article~7.---14~P.}
\bib{Karp~D., Sitnik~S.~M.}{Asymptotic approximations for the first incomplete elliptic integral near logarithmic singularity~/\!/~2006.---Arxiv:~math.~CA/0604026.---20~P.}
\bib{Karp~D., Savenkova~A., Sitnik~S.~M.}{Series expansions and asymptotics for incomplete elliptic integrals
via partial fraction decompositions~/\!/~Proceedings of the fifth annual conference of the
Society for special functions and their applications (SSFA).---2004.---Lucknow (India).---P.~4--30.}
\bib{Kazarinoff~N.~D.}{Analytic Inequalities.---(1st ed. 1961), Holt, 2003---89~p.}
\bib{Kufner~A., Maligranda~L., Persson~L.-E.}{The Hardy Inequality.---Pilsen,~2007.---162~p.}
\bib {Kusraev~A.~G., Buskes~G.}{Representation and extension
of orthoregular bilinear operators~/\!/Владикавказский математический журнал.---2007.---Т.~9, \No~1.---С.~1--17.}
\bib {Kusraev~A.~G.}{H\"{o}lder type inequalities for
orthosymmetric bilinear operators~/\!/Владикавказский математический журнал.---2007.---Т.~9, \No~3.---С.~3--37.}

\bib {Maligranda~L.}{Why Holder's inequality should be called Roger's inequality~/\!/Math. Inequal. Appl.---1998.---\No~1.---P.~69--83.}
\bib{Mikhlin~S.~G., Lehmann~R.}{Constants in Some Inequalities of Analysis.---Wiley, 1986.---108~p.}
\bib{Milovanovi\'{c}~G.~V.}{Recent Progress in Inequalities. (Volume is dedicated to Professor Dragoslav S. Mitrinovi\'{c} (1908--1995).---Kluwer, 1998.---532~p.}
\bib{Mitrinovi\'{c}~D.~S.}{Elementary Inequalities.---Groningen, 1964.---159~p.}
\bib{Mitrinovi\'{c}~D.~S.}{Nejednakosti.---Beograd, 1965.---240~p.}
\bib{Mitrinovi\'{c}~D.~S., Vasi\'{c}~P.~M.}{Sredine.---Beograd, 1969.---122~p.}
\bib{Mitrinovi\'{c}~D.~S. (in cooperation with P.M.\,Vasi\'{c}).}{Analytic Inequalities.---Springer, 1970.---400~p.}
\bib{Mitrinovi\'{c}~D.~S., Pe\v{c}ari\'{c}~J.~E., Volenec~V.}{Recent Advances in Geometric Inequalities.---Springer, 1989.---736~p.}
\bib{Mitrinovi\'{c}~D.~S., Pe\v{c}ari\'{c}~J.~E., Fink~A.~M.}{Classical and new inequalities in analysis.---Kluwer,
1993.---740~p.}
\bib{Niculescu~C., Persson~L.~E.}{Convex functions and their applications: a contemporary approach.---Springer,~2006.---255~p.}

\bib{Opic~B., Kufner~A.}{Hardy--Type Inequalities.---Longman,1990.---333~p.}

\bib{Pales~Zs.}{Inequalities for differences of powers~/\!/Math. Anal. Appl.---1988.---\No~131.---P.~271--281.}
\bib{Pales~Zs.}{Comparison of two variable homogeneous means~/\!/General
Inequalities 6. Proc. 6th Internat. Conf. Math. Res. Inst.
Oberwolfach.---Basel: Birkhauser Verlag, 1992.---P.~59--69.}
\bib{Rado~T.}{On convex functions~/\!/Trans. Amer. Math. Soc.---1935.---Vol.~37.---P.~266--285.}
\bib{S\'{a}ndor~J.}{Sеlected chapters of geometry, analysis and
number theory.---RGMIA monographs, Victoria University, 2005.---460~p.}
\bib{Sitnik~S.~M.}{Refinements of the Bunyakovskii -- Schwartz inequalities with applications to special functions estimates~/\!/~ Conference in Mathematical Analysis and Applications in Honour of Lars Inge Hedberg's 60 th Birthday. Link\"{o}ping University, Sweden, 1996.---P.~97.}
\bib{Steele~J.~M.}{The Cauchy-Schwarz Master Class: An Introduction to the Art of Mathematical Inequalities.---Cambridge University Press, 2004.---306~p.}
\bib{Taneja~I.~J.}{Generalized Information Measures and Their Applications.---2001.---on--line book: www.mtm.ufsc.br/~taneja/book/book.html}
\bib{Toader~Gh., Toader~S.}{Greek Means and the Arithmetic-Geometric Mean.---RGMIA Monographs, Victoria University, 2005.---95~p.}

\bib{Zarantonello~E.~H.}{The meaning of the Cauchy-Schwarz-Buniakovsky inequality (MRC Technical summary report).}

\end{enumerate}

 \Adress{\textsc{Ситник Сергей Михайлович}\\
 Воронежский институт МВД России;\\
 кафедра высшей математики.\\
 Россия, 394065, Воронеж, Проспект патриотов, 53\\
 E-mail: \verb"box2008in@gmail.com"}

 \begin{center}
 Sitnik S. M.

 Generalizations of classic inequalities
 \end{center}

\Abstract{In this survey we consider generalizations and specifications of classic inequalities. First  generalizations to the Young inequality is studied with possible specifications and unsolved problems. The main part consists of authour's method to generalize Cauchy--Bunyakowsky  inequalities both in discrete and integral cases. Many applications are also considered.}

\end{document}